\documentclass[review,onefignum,onetabnum]{siamonline220329}

\usepackage{lipsum}
\usepackage{amsfonts}
\usepackage{graphicx}
\usepackage{epstopdf}
\usepackage{algorithmic}
\ifpdf
  \DeclareGraphicsExtensions{.eps,.pdf,.png,.jpg}
\else
  \DeclareGraphicsExtensions{.eps}
\fi

\usepackage{enumitem}
\setlist[enumerate]{leftmargin=.5in}
\setlist[itemize]{leftmargin=.5in}

\newsiamremark{remark}{Remark}
\newsiamremark{hypothesis}{Hypothesis}
\crefname{hypothesis}{Hypothesis}{Hypotheses}
\newsiamthm{claim}{Claim}

\headers{Boltzmann equation with uncertainties}{L.~Liu and K.~Qi}

\title{Spectral convergence of a semi-discretized numerical system for the spatially homogeneous Boltzmann equation with uncertainties \thanks{Submitted to the editors DATE.
}}

\author{Liu Liu\thanks{Department of Mathematics, The Chinese University of Hong Kong, Shatin, Hong Kong, People's Republic of China.
  (\email{lliu@math.cuhk.edu.hk}).}
\and Kunlun Qi\thanks{School of Mathematics, University of Minnesota--Twin Cities, Minneapolis, MN 55455 USA.
  (\email{kqi@umn.edu}).}
}

\usepackage{amsopn}

\usepackage{amsmath,amssymb}

\usepackage{bm}
\usepackage{cite}
\usepackage{comment}
\usepackage[page,toc,titletoc]{appendix}
\usepackage{todonotes}

\usepackage{bm}

\newcommand*{\im}{\mathop{}\!\mathrm{i}}
\newcommand*{\e}{\mathop{}\!\mathrm{e}}

\newcommand{\rd}{\mathrm{d}}

\newcommand{\mP}{\mathcal{P}}

\newcommand{\bR}{\mathbb{R}}
\newcommand{\bS}{\mathbb{S}}

\newcommand{\Sd}{\mathbb{S}}
\newcommand{\D}{\mathcal{D}_L}

\newcommand{\fe}{f_{N}}

\newcommand{\n}{\bm{\mathrm{n}}}
\newcommand{\kk}{\bm{\mathrm{k}}}
\newcommand{\lb}{\bm{\mathrm{l}}}
\newcommand{\m}{\bm{\mathrm{m}}}
\newcommand{\ib}{\bm{\mathrm{i}}}
\newcommand{\jb}{\bm{\mathrm{j}}}

\ifpdf
\hypersetup{
  pdftitle={Spectral convergence of a semi-discretized numerical system for the spatially homogeneous Boltzmann equation with uncertainties},
  pdfauthor={L.~Liu and K.~Qi}
}
\fi

\begin{document}
	\nolinenumbers

\maketitle

\begin{abstract}
In this paper, we study the Boltzmann equation with uncertainties and prove that the spectral convergence of the semi-discretized numerical system holds in a combined velocity and random space, where the Fourier-spectral method is applied for approximation in the velocity space whereas the generalized polynomial chaos (gPC)-based stochastic Galerkin (SG) method is employed to discretize the random variable. Our proof is based on a delicate energy estimate for showing the well-posedness of the numerical solution as well as a rigorous control of its negative part in our well-designed functional space that involves high-order derivatives of both the velocity and random variables. This paper rigorously justifies the statement proposed in [Remark 4.4, J.~Hu and S.~Jin, \textit{J.~Comput.~Phys.},~315 (2016),~pp.~150-168].
\end{abstract} 





\section{Introduction}
\label{sec:intro}

\textbf{Background and goals} The study of kinetic equations has undergone a long history, due to their wide applications in diverse and crucial scientific domains, encompassing rarefied gas dynamics, plasma physics, astrophysics, and emerging fields like semiconductor device modeling \cite{Markowich}. Additionally, they find application in environmental, social, and biological sciences \cite{BGK-book}, providing a comprehensive framework to describe the non-equilibrium dynamics of gases or systems comprising a multitude of particles. 
As a typical example, the Boltzmann equation is the commonly used governing equation in modeling different phenomena, where it is adept at capturing rarefied gas flows pertinent to hypersonic aerodynamics, gases in vacuum technologies, and fluid dynamics within microelectromechanical devices \cite{cercignani}. Recently, the Boltzmann equation also extended its utility to the modeling of social and biological phenomena \cite{PT2, AHP2015, KW2000}. For comprehensive insights and up-to-date references, we refer the readers to \cite{DPR,pareschi2001introduction, rjasanow, Villani02}.

Albeit the significance in real applications, the numerical approximation of the Boltzmann equation has been an extremely challenging problem due to its high dimension, non-linearity, and non-locality. 
To overcome the curse-of-dimensionality, the particle-based direct simulation Monte Carlo (DSMC) method \cite{Nanbu80,Bird} has demonstrated notable advantages in terms of efficiency and ease of implementation, however, this type of stochastic method usually suffered from the slow convergence as well as the inefficiency particularly when applying to the non-steady and low-speed flow. 
In contrast, deterministic type of methods have experienced rapid development in recent decades, fueled by advancements in computing power \cite{DP15}. Notably, for approximating the Boltzmann collision operator, the Fourier-Galerkin spectral method offers a compelling framework \cite{PP96, PR05}. This method not only achieves spectral accuracy but also lends itself to easy adaptation through accelerated algorithms \cite{MP06, MPR13, GHHH17}.
Beyond the accomplishments of the Fourier-Galerkin spectral method in numerical simulation, the attainment of a rigorous theoretical proof regarding its spectral convergence has remained elusive for an extended period, necessitating the incorporation of filters to maintain solution positivity \cite{PR00stability} Inspired by the pioneering stability analysis in \cite{FM11}, which leverages the ``spreading" property of the collision operator, the authors in recent work \cite{HQY21} introduce a novel framework to demonstrate the stability and convergence of the Fourier-Galerkin spectral method. The approach relies on a meticulous $L_v^2$-estimate of the negative component of solutions to the spatially homogeneous Boltzmann equation. Besides, it is worthwhile mentioning the recent advancement in developing the conservative spectral method \cite{AGT18, PR2022, CLL2024}, which naturally enables the implication of spectral convergence and stability.

On the other hand, there are usually various sources of uncertainties in the kinetic equations such as modeling errors, imprecise measurements, and uncertain initial conditions. 
Studying the uncertainty quantification (UQ) problem is crucial to assess, validate and enhance the underlying models, which makes our project imperative. 
In particular, the collision kernel or scattering cross-section in the Boltzmann equation delineates the transition rate occurring during particle collisions. The task of calculating this collision kernel from first principles is exceptionally intricate. In practice, only heuristic approximations and empirical data are used, which inevitably introduce uncertainties into the collision kernel. Furthermore, uncertainties may arise from inaccurate measurements of initial or boundary data, as well as from source terms, bringing additional uncertainties to the model. To numerically study the Boltzmann equation and other kinetic models with randomness, we refer the readers to \cite{HPY, JXZ15, Jin-ICM,  Poette2, RHS, LW2017, CMT18, PZ2020}. Among different numerical algorithms, the generalized polynomial chaos (gPC)-based stochastic Galerkin (SG) method and its variations have been adapted widely and shown success in broad applications \cite{Xiu}.
Aside from numerical simulation, we mention the relevant works for the theoretical proof of the stability and convergence in random approximation: the spectral convergence has been proved for the proposed gPC-based SG method in \cite{HJ16}, while in \cite{LJ18} and \cite{DJL19}, the authors established powerful framework based on hypocoercivity to conduct the local sensitivity analysis for a class of multiscale, inhomogeneous kinetic equations with random uncertainties, the approximation of which is utilized by the gPC-based SG method. We refer to the recent collection \cite{JinPareschi}, survey \cite{parUQ} and some relevant works \cite{PIN_book, LeMK, Bio-UQ, NTW, PDL, Schwab}. \\[-6pt]

\textbf{Challenges and our contributions} 
Although the convergence and stability of the numerical solution to the Boltzmann equation has been extensively shown in either the deterministic case \cite{PR00stability, FM11, HQY21, CLL2024} by the Fourier spectral method, or the stochastic case via gPC-based stochastic Galerkin method \cite{HJ16, Jin-ICM, GHY2021}, there are very few results concerning the semi-discretized system that simultaneously takes the deterministic and stochastic approximation into consideration, due to the complex structure of the nonlinear Boltzmann operator as well as the lack of the suitable functional space to work in. Given the motivations above, in this paper, our main purpose is to show that spectral convergence and stability still hold for the numerical approximation of the Boltzmann equation with uncertainties in the combined velocity and random space, which can be seen as a rigorous justification of the statement proposed by Hu-Jin in \cite[Remark 4.4]{HJ16}.
Specifically speaking, to show the convergence and stability, i.e., quantify the ``error" in an appropriate functional space (usually the Sobolev space), a key point is to study the well-posedness of the numerical solution obtained from the semi-discretized system \eqref{PFSNK} in the same class of space; to achieve this, it will heavily rely on the delicate estimate of the collision operator in the higher-order Sobolev space, however, some well-known $L^p$- and $H^k$-estimates of the collision operator without uncertainties cannot be directly applied due to the existence of the random variable. Besides, the intrinsic drawback of the Fourier spectral method, e.g., some important physical properties such as positivity are unable to be conserved, remains the main trouble-maker in deriving the well-posedness of the numerical solution, which is even harder than the merely deterministic case. 

To overcome the aforementioned difficulties, following up our previous work \cite{LQ2022}, we first present our brand new estimates of the collision operator in the Sobolev space that involves both the velocity and random variables (see Section \ref{subsubsec:QR}); since only the uniform bound of collision kernel is assumed in the random variable, we have to complete the estimate in the high-order Sobolev space by considering the embedding relation. In terms of the lack of conservative quantities, namely, the positivity, we apply a careful Sobolev estimate of the negative part of the numerical solution to make sure it is under control at least within a prescribed time interval.
Such a strategy has proved to work well in the deterministic case \cite{HQY21}, however, we need to show that it still is valid in our new functional space (see Section \ref{subsec:propagation}). \\[-6pt]

\textbf{Organization of our paper} The rest of this paper is organized as follows: In Section \ref{sec:Boltzmann}, we first introduce the formulation of the semi-discretized numerical system. Then, some preliminary assumptions of the collision kernel and initial condition are presented in Section \ref{sec:pre}, along with some useful estimates and properties that will be utilized throughout the paper. Furthermore, we show the well-posedness of the numerical solution with uncertainties in our designated functional space in Section \ref{sec:main}. Based on the well-posedness, the spectral convergence of the obtained numerical solution will be proved in Section \ref{sec:conv}.


\section{The Boltzmann equation and semi-discretized numerical system}
\label{sec:Boltzmann}

\subsection{The Boltzmann equation with uncertainties}
\label{subsec:Boltzmann_uncertain}

We consider the spatially homogeneous Boltzmann equation with uncertainties, 
\begin{equation}\label{IBE}
\partial_{t} f (t,v,z) = Q(f,f)(t,v,z), \quad t>0, \ v\in \mathbb{R}^d,\, d\geq 2, \quad z \in I_{z} \subset \mathbb{R}^{d_z},
\end{equation}
with the initial datum
\begin{equation}\label{f0}
f(0,v,z) = f^{0}(v,z),
\end{equation}
where $f=f(t,v,z)$ is the probability density function at time $t$, with velocity variable $v$ and $d_z$-dimensional random variable $z$ that is assumed known and characterizes the random inputs. 
Here $Q$ is the bilinear collision operator describing binary collisions among particles, 
\begin{equation} \label{Qstrong}
Q(g,f)(v, z)=\int_{\bR^d}\int_{\Sd^{d-1}}B(|v-v_*|,\cos \theta, z)[g(v_*',z)f(v',z)-g(v_*,z)f(v,z)] \,\rd{\sigma}\, \rd{v_*}, 
\end{equation}
where $\sigma$ is a unit vector varying over the sphere $\Sd^{d-1}$ such that the post-collisional velocities $v'$, $v_*'$ are defined as
\begin{equation}\label{v'vs'}
v'=\frac{v+v_*}{2}+\frac{|v-v_*|}{2}\sigma, \quad v_*'=\frac{v+v_*}{2}-\frac{|v-v_*|}{2}\sigma.
\end{equation}
As indicated above, in this paper, we consider uncertainties coming from 
\begin{itemize}
\item[(i)] the initial datum $ f^{0}(v,z)$;
\item[(ii)] the collision kernel $B(|v-v_*|,\cos\theta, z) \geq 0$ that owns the form: 
\begin{equation}\label{kernel}
B(|v-v_*|,\cos\theta, z) = \Phi(|v-v_*|)\ b(\cos\theta, z), \qquad \cos\theta =\frac{\sigma\cdot (v-v_*)}{|v-v_*|},
\end{equation}
where the kinetic part $\Phi$ is a non-negative function, the angular part $b$ satisfies the Grad's cut-off assumption, that is, for all $z \in I_{z} \subset \mathbb{R}^{d_z}$, 
\begin{equation}\label{cutoff}
\int_{\Sd^{d-1}} b(\cos\theta, z) \,\rd{\sigma} < \infty.
\end{equation}
One typical collision kernel that is widely used is the variable hard sphere (VHS) model \cite{HJ16}:
\begin{equation}\label{VHS}
    B(|v-v_*|, \cos\theta, z) = b(z) |v-v_*|^{\gamma}, \quad 0 \leq \gamma \leq 1.
\end{equation}
\end{itemize}
The random variable $z$ comprises random vectors in $d_z$ dimension, each adhering to a known distribution. To simplify our analysis, we assume that all components of $z$ are mutually independent and have been acquired through a dimension reduction technique, such as the Karhunen-Loeve expansion \cite{Loeve1977, Xiu}:
\begin{equation}
    b(\cos\theta,z) \approx b_0(\cos\theta) + \sum_{i=1}^{d_z} b_i(\cos\theta) z_i,
\end{equation}
with $z_1, \cdots, z_{d_z}$ independent random variables following the probability density function $\pi(z)$.

\subsection{Formulation of the semi-discretized numerical system}
\label{sec:review}

We introduce the formulation of our semi-discretized (except for the discretization in the temporal space) numerical system for solving the spatially homogeneous Boltzmann equation with uncertainties \eqref{IBE}, where both deterministic and stochastic approximations are imposed simultaneously. More precisely, we apply a Fourier-Galerkin spectral method for velocity variable $v$ \cite{PR00, DP15, MP06, GHHH17, HQ20}, while a gPC-based stochastic Galerkin method is utilized in the discretization of random variable $z$ \cite{HJ16, LJ18, Xiu, XK02}.

To apply the Fourier-Galerkin spectral method in velocity space, we consider an approximation of \eqref{IBE} on a periodic domain of velocity $\mathcal{D}_L=[-L, L]^d$,
\begin{equation}\label{ABE}
\left\{
\begin{aligned}
&\partial_{t} f (t,v,z) = Q^{R}(f,f)(t,v,z), \quad t>0, \ v\in \mathcal{D}_L, \, z \in I_{z} \subset \mathbb{R}^{d_z},\\[4pt]
& f(0,v,z)= f^{0}(v,z), 
\end{aligned}
\right.
\end{equation}
where the initial condition $f^0$ is a non-negative periodic function, $Q^{R}$ is the truncated collision operator of \eqref{Qstrong},
\begin{equation}\label{QR}
\begin{split}
Q^R(g,f)(v,z)&:=\int_{\mathcal{B}_R}\int_{\Sd^{d-1}}\Phi(|q|)\, b(\sigma\cdot \hat{q},\,z)\left[g(v'_*,z)f(v',z)-g(v-q,z)f(v,z)\right] \rd{\sigma}\, \rd{q}\\[4pt]
&=\int_{\bR^d} \int_{\Sd^{d-1}} B_T(|q|,\sigma\cdot\hat{q},z) \left[g(v'_*,z)f(v',z)-g(v-q,z)f(v,z)\right]\rd{\sigma}\, \rd{q},
\end{split}
\end{equation}
where, after the change of variable $v_* \mapsto q=v-v_*$, the relative velocity $q$ is truncated to a ball $\mathcal{B}_R$ with radius $R$ centered at the origin; $B_T(|q|,\sigma\cdot\hat{q},z)$ is the truncated collision kernel $B_T(|q|,\sigma\cdot\hat{q},z):= \mathbf{1}_{|q|\leq R}\,\Phi(|q|)\,b(\sigma\cdot \hat{q}, z)$ by denoting $q=|q|\hat{q}$ with $|q|$ the magnitude and $\hat{q}$ the directional unit vector.
The values of $L$ and $R$ are positive satisfying $L \geq R > 0$,
where, to avoid aliasing errors \cite{PR00}, one usually takes
\begin{equation}\label{RL}
R=2S, \quad L\geq \frac{3+\sqrt{2}}{2}S,
\end{equation}
where the support set of $f^0(v,z)$ in $v$ is assumed to be within $\mathcal{B}_S$ for all $z \in I_z$.

Now we are in a position to present the semi-discretized numerical system, i.e., we seek an approximated solution in the following form: 
\begin{equation}\label{fNK}
	f(t,v,z) \approx f_{N}^{K}(t,v,z) : = \sum_{|\kk|= 0}^{K} \sum_{|\n|= 0}^{N} f^{\kk}_{\n}(t)\Phi_{\n}(v) \Psi^{\kk}(z),
\end{equation}
where $N$ and $K$ are non-negative integers, $\n = (n_1,\dots,n_d)$ and $ \kk = (k_{1},...,k_{d_{z}}) $ are multi-indexes with $|\n| = n_1+\dots+n_d$ and $ |\kk| = k_{1}+...+k_{d_{z}}$, and
\begin{equation*}
\Phi_{\n}(v) \in \mathbb{P}_N :=\text{span} \left\{  \e^{\im \frac{\pi}{L} \n \cdot v} \, \Big| \, 0 \leq |\n| \leq N \right\} 
\end{equation*}
The space $\mathbb{P}_N$ is equipped with the inner product in $v$: 
\begin{equation*}
\langle f(t,\cdot,z), g(t,\cdot,z) \rangle_{\D} = \frac{1}{(2L)^{d}}\int_{\D}f(t,v,z) \, \overline{g(t,v,z)} \,\rd v.
\end{equation*}
while $\Psi^{\kk}(z)$ are the orthogonal gPC basis functions satisfying
\begin{equation*}
\int_{I_{z}} \Psi^{\ib}(z) \, \Psi^{\jb}(z)  \, \pi(z) \,\rd z = \delta_{\ib\jb}, \quad 0 \leq |\ib|, |\jb| \leq K 
\end{equation*}
with $\pi(z)$ being the the probability distribution function of $z$. The space 
$$\mathbb{P}^K := \text{Span} \left\{ \Psi^{\kk}(z) \, \Big|\, 0 \leq |\kk| \leq K \right\}$$
is equipped with the inner product in $z$:
\begin{equation*}
	\langle f(t,v,\cdot), \, g(t,v,\cdot) \rangle_{I_{z}} = \int_{I_{z}}  f(t,v,z) \, \overline{g(t,v,z)} \, \pi(z) \,\rd z.
\end{equation*}

Inserting \eqref{fNK} into \eqref{IBE} and conducting projections onto the space $\mathbb{P}_N$ and $\mathbb{P}^K$ successively yields
\begin{equation} \label{PFSNK}
	\left\{
	\begin{aligned}
		&\partial_{t} f_{N}^{K}(t,v,z) = \mP_N^K Q^{R}(f_N^K,f_N^K)(t,v,z), \quad t>0, \ v\in \mathcal{D}_L,\ z \in I_{z}, \\[4pt]
		& f_N^K(0,v,z)=f_{N}^{0,K}(v,z),
	\end{aligned}
	\right.
\end{equation}
where $\mP_N^K$ is the projection on any suitable function $g$ such that
\begin{equation}\label{projNK}
	\mP_N^K g(t,v,z) = \sum_{|\kk|= 0}^{K} \sum_{|\n|= 0}^{N}  g^{\kk}_{\n}(t) \Phi_{\n}(v) \Psi^{\kk}(z), \quad g_{\n}^{\kk}(t) = \Big \langle \langle g(t,\cdot,\cdot), \Phi_{\n}(\cdot) \rangle_{\D}, \Psi^{\kk}(\cdot) \Big\rangle_{I_{z}}. 
\end{equation}
More discussions on the initial approximation $f_{N}^{0,K}$ will be given in Section \ref{subsec:initial}. We also refer the readers for a detailed introduction of the separate projection $\mP_N$ and $\mP^K$ in our preceding work \cite[Section 2.2]{LQ2022}.

Furthermore, we present the semi-discretized system in its dual mode: 
\begin{equation} \label{FSNK}
	\left\{
	\begin{aligned}
		&\partial_{t} f^{\kk}_{\n}(t) = Q^{R,\kk}_{\n}(f^{K}_{N},f^{K}_{N})(t), \quad 0 \leq |\n| \leq  N, \quad  0 \leq |\kk| \leq  K, \\[4pt]
		& f^{\kk}_{\n}(t=0,v,z)= f^{0,\kk}_{\n},
	\end{aligned}
	\right.
\end{equation}
with
\begin{equation*}
\begin{aligned}
 & Q^{R,\kk}_{\n}(f_{N}^{K},f_{N}^{K}):= \Big \langle \langle Q^{R}(f^{K}_{N},f^{K}_{N})(t,\cdot,\cdot), \Phi_{\n}(\cdot) \rangle_{\D}, \Psi^{\kk}(\cdot) \Big\rangle_{I_{z}}, \\[2pt]
& f^{0,\kk}_{\n}:= \Big \langle \langle f^K_{N}(0,\cdot,\cdot), \Phi_{\n}(\cdot) \rangle_{\D}, \Psi^{\kk}(\cdot) \Big\rangle_{I_{z}}.
\end{aligned}
\end{equation*}
In particular, we can deduce the term $Q^{R,\kk}_{\n}$ as follows:
\begin{equation}\label{QRnk}
	Q^{R,\kk}_{\n}(t) = \sum_{|\ib|,|\jb| = 0}^{K} S^{\kk\ib\jb} \sum\limits_{\substack{|\lb|,|\m|=0 \\[2pt] \lb+\m=\n}}^{N} G(\lb,\m)\, f^{\ib}_{\lb}\, f^{\jb}_{\m}
\end{equation}
where the weight $G(\lb,\m)$ and $S^{\kk\ib\jb}$ are given by
\begin{equation}\label{GG}
\begin{split}
G(\lb,\m) & = \int_{\mathcal{B}_{R}}\int_{\Sd^{d-1}}\Phi(|q|)b(\sigma\cdot \hat{q}) \left[ \e^{-\im \frac{\pi}{2L}(\lb+\m)\cdot q +\im \frac{\pi}{2L}|q|(\lb-\m)\cdot \sigma} - \e^{-\im \frac{\pi}{L}\m\cdot q} \right] \rd\sigma\,\rd q\\[4pt]
& =\int_{\mathcal{B}_{R}}\e^{-\im \frac{\pi}{L}\m \cdot q}\left[\int_{\Sd^{d-1}}\Phi(|q|)b(\sigma\cdot \hat{q})(\e^{\im \frac{\pi}{2L}(\lb+\m)\cdot (q-|q|\sigma)}-1)\, \rd\sigma\right]\rd q, 
\end{split}
\end{equation}
and
\begin{equation}\label{S_kij}
S^{\kk\ib\jb} := \int_{I_{z}} b(z) \Psi^{\kk}(z) \Psi^{\ib}(z) \Psi^{\jb}(z) \pi(z) \,\rd z.
\end{equation}
It is crucial to emphasize that the evaluations of $G(\lb,\m)$ and $S^{\kk\ib\jb}$ are precomputed processes, undertaken only once for a given collision kernel. In addition, for a specific category of collision kernels, such as the VHS model, both $G(\lb,\m)$ and $S^{\kk\ib\jb}$ can be simplified to more concise forms \cite{GHHH17, HJ16}, facilitating easier evaluation.



For more details and related fast algorithms about applying the Fourier Spectral method and gPC-based stochastic Galerkin method for the Boltzmann equation with uncertainties, we refer the readers to \cite{DP15, MP06, GHHH17, HQ20, HJ16, LJ18} and the references therein.



\section{Preliminary}
\label{sec:pre}

\subsection{Notations}
\label{subsec:norms}
For a function $f(t, v, z)$ that is periodic in velocity space $\D$, we define its Lebesgue norm and Sobolev norm w.r.t the velocity variable:
\begin{equation}
\|f(t,\cdot,z)\|^p_{L^p_{v}}:= \int_{\D} |f(t, v, z)|^p \,\rd{v}, \quad 
\|f(t, \cdot, z)\|^2_{{H^k_{v}}}:= \sum_{|\nu|\leq k} \|\partial_v^{\nu} f(t,\cdot,z) \|^2_{L^2_{v}}, 
\end{equation}
and the norms w.r.t the random variable in $I_z$:
\begin{equation}
\|f(t,v,\cdot)\|^p_{L^p_{z}}:= \int_{I_z} |f(t, v, z)|^p \, \pi(z) \,\rd z, 
\quad 
\|f(t, v, \cdot)\|^2_{{H^r_{z}}}:= \sum_{|\mu|\leq r} \|\partial_z^{\mu} f(t,v,\cdot) \|^2_{L^2_{z}}, 
\end{equation}
as well as the high-order Sobolev norms combining the velocity and random variables:
\begin{equation}
    \| f(t,\cdot,\cdot) \|^2_{H^k_v H^r_z} := \sum_{|\nu|\leq k} \sum_{|\mu|\leq r} \| \partial_v^\nu \partial_z^\mu f(t,\cdot,\cdot)\|^2_{L_v^2 L_z^2}.
\end{equation}

Similar to the deterministic case, for a function $f(v, z)$ and each $z \in I_z$, 
the positive and negative parts are defined by 
\begin{equation}
f^+(v, z)=\max\limits_{v\in\D}\{ f(v, z), 0 \}, \quad 
f^-(v, z)=\max\limits_{v\in\D}\{ -f(v, z),0 \}.
\end{equation}

\subsection{Assumptions on the collision kernel and initial condition}
\label{subsec:initial}

\noindent{\bf Basic assumptions on the collision kernel:} 

\begin{itemize}

\item[(i)] Considering the collision kernel in the form \eqref{kernel}, the kinetic part $\Phi$ satisfies
\begin{equation}\label{kinetic}
\left\| \mathbf{1}_{|q|\leq R}\Phi(|q|)\right\|_{L^{\infty}(\D)} < \infty . 
\end{equation}
Notice that the so-called hard potentials $\Phi(|q|)=|q|^{\gamma}$ ($\gamma \geq 0$) and the ``modified" soft potentials $\Phi(|q|)=(1+|q|)^{\gamma}$ ($-d<\gamma <0$) all satisfy this condition. 

\item[(ii)] For all $z\in I_z$, we assume a uniform bound for all the $z$-derivatives of $b$: 
\begin{equation}\label{Assump_B}
 b(\sigma\cdot\hat{q},\,z) > 0, \qquad
 |\partial_z^l b(\sigma\cdot\hat{q},\,z)| \leq C_b \quad \text{for} \quad
 0\leq |l| \leq r. 
\end{equation}
\end{itemize}

\begin{remark}
We mainly focus on the collision kernel with uncertainties in the angular kernel $b(\cos\theta, z)$ in the work. 
While beyond the scope of this paper, it is noteworthy that our analysis can be similarly extended to the case where the kinetic part $\Phi$ in the collision kernel is assumed uncertain, i.e.,
$ B_T(|q|,\,\sigma\cdot\hat{q},\,z) = \mathbf{1}_{|q|\leq R} \Phi(|q|,\,z)\, b(\sigma\cdot\hat{q})$, if more stringent conditions are carefully imposed on $\Phi(|q|,\,z)$ to ensure the uniform boundedness of $\partial_z^l B_T$ for all $z \in I_z$ and $0 \leq |l| \leq r$.

\end{remark}
\bigskip

\noindent{\bf Basic assumptions on the initial condition:} 

To prove the well-posedness and stability of the numerical solution $f_N^K$ to \eqref{PFSNK}, we need to restrict to a certain class of initial data. For the initial condition $f^0(v,z)$ to the continuous problem \eqref{IBE}, we assume it to be non-negative, periodic in the velocity space, and belongs to $L_v^1 H_z^1 \cap H_v^1 H_z^1$. For the approximated initial condition $f_N^{0,K}(v,z) = \mP_N^K f^0(v,z) $ \footnote{Note that $f_N^{0,K}=\mP_N^K f^0$ is not the only possible initial condition for the numerical system, any reasonable numerical approximation satisfying the following assumptions (i)--(iv) will work.} to the semi-discretized numerical system \eqref{PFSNK}, we can show that it satisfies the following properties \cite{HQY21, Xiu}: 
\begin{itemize}
\item[(i)] Mass conservation of the approximation: we assume, for $r \geq 0$,
\begin{equation} \label{con(a)}
\| \int_{\D} f^{K,0}(v,\cdot) \,\rd v \|_{H^r_z} = \| \int_{\D} f^0(v,\cdot)\, \rd v \|_{H^r_z}.
\end{equation}
The conservation of $\| \int_{\D} f_N^K(t,v,\cdot)  \,\rd v \|_{H_z^r}$ can be verified by further considering the Lemma \ref{lemma:conv}.

\item[(ii)] Control of $L_v^2 H^1_z$ and $H_v^1 H^1_z$ norms: for any integer $N, K \geq 0$,
\begin{equation} \label{con(b)}
\|f^{0,K}_N(\cdot,\cdot)\|_{L_v^2 H^1_z} \leq \|f^0(\cdot,\cdot)\|_{L_v^2 H^1_z}, \quad \|f^{0,K}_N(\cdot,\cdot)\|_{H_v^1 H^1_z}\leq \|f^0(\cdot,\cdot)\|_{H_v^1 H^1_z}.
\end{equation}
In fact, for $f^0 \in H_v^1 H_z^1 $, the projected initial condition $f^{0,K}_N = \mathcal{P}_N^K f^0$ automatically satisfies this condition with $\mathcal{P}_N^K$ defined in \eqref{projNK}, which is a direct consequence of the Parseval type identity.

\item[(iii)] Control of $L_v^1 H^1_z$ norm: there exists integers $N_0$ and $K_0$ such that for all $N > N_0$ and $K > K_0$,
\begin{equation} \label{con(c)}
\|f_N^{0,K}(\cdot,\cdot)\|_{L_v^1 H^1_z} \leq C \|f^0(\cdot,\cdot)\|_{L_v^1 H^1_z}.
\end{equation}

\item[(iv)] $L_v^2 H^1_z$-norm of $f_N^{0,K,-}$ can be approximated arbitrarily small: for any $\varepsilon>0$, there exists integers $N_0$ and $K_0$ such that for all $N > N_0$ and $K > K_0$,
\begin{equation} \label{con(d)}
\|f_N^{0,K,-}(\cdot,\cdot)\|_{L_v^2 H^1_z} < \varepsilon,
\end{equation}
where $f_N^{0,K,-}$ stands for the negative part of initial datum $f_N^{0,K}$. This condition can be guaranteed, if the uniform convergence of the Fourier expansion and the gPC expansion are satisfied; to achieve it, one might impose continuity on $f^0$. For example, it would be sufficient if $f^0$ is H\"older continuous, or continuous plus bounded variation in the velocity and random space (In fact, thanks to the embedding relation, the continuity requirement will be satisfied, since $f^0$ lies in the higher-order Sobolev space in our convergence analysis).
\end{itemize}


\subsection{Useful properties}
\label{subsec:pre}

In this subsection, we will present some useful properties that will be repeatedly applied in the following proofs.

\subsubsection{Estimates for the collision operator \texorpdfstring{$Q^{R}$}{QR}}
\label{subsubsec:QR}

The estimates of the truncated collision operator $Q^R$ with uncertainty will play a key role in the proof for our main results. Under the cut-off assumption \eqref{cutoff}, it is more convenient to show the estimates of $Q^R$ by splitting it into the gain term $Q^{R,+}$ and loss term $Q^{R,-}$, which are defined as follows:
\begin{equation}
\begin{split}
 Q^{R,+}(g,f)(v,z): &=\int_{\bR^d}\int_{\Sd^{d-1}}\mathbf{1}_{|q|\leq R}\,\Phi(|q|) \,b(\sigma\cdot \hat{q},z)\, g(v'_*,z)f(v',z) \,\rd{\sigma} \,\rd{q}, \\[4pt]
 Q^{R,-}(g,f)(v,z): &=\int_{\bR^d}\int_{\Sd^{d-1}}\mathbf{1}_{|q|\leq R}\,\Phi(|q|)\,b(\sigma\cdot \hat{q},z)\,g(v-q,z)f(v,z) \,\rd{\sigma} \,\rd{q}\\[4pt]
 & = f(v,z)L^R[g](v,z),
\end{split}
\end{equation}
where
\begin{equation*}
    L^R[g](v,z):= \int_{\bR^d}\int_{\Sd^{d-1}}\mathbf{1}_{|q|\leq R}\,\Phi(|q|)\,b(\sigma\cdot \hat{q},z)\,g(v-q,z) \,\rd{\sigma} \,\rd{q}.
\end{equation*}

Now we first present the $L^p_{v,z}$ estimate of $Q^R$ in the Lebesgue space that encompasses both velocity and random variables. 

\begin{proposition}\label{pro_QLp}
    Let the collision kernel $ B $ satisfy the assumptions \eqref{kernel}-\eqref{cutoff}, \eqref{kinetic}-\eqref{Assump_B} and truncation parameters $R$ and $L$ satisfy \eqref{RL}, then the truncated collision operators $ Q^{R,\pm}(g,f)$ satisfy the following estimates: for $1\leq p \leq \infty$,	
    \begin{equation}\label{QGLp1}
    \begin{split}
     \left\|Q^{R,+}(g,f)(\cdot,\cdot)\right\|_{L^p_{v,z}} \leq C^+_{R,L,I_z,d,p}(B) \left\|g(\cdot,\cdot)\right\|_{L^{1}_v L^{\infty}_z} \left\|f(\cdot,\cdot) \right\|_{L^p_{v,z}},\\[4pt]
     \left\|Q^{R,+}(g,f)(\cdot,\cdot)\right\|_{L^p_{v,z}} \leq C^+_{R,L,I_z,d,p}(B) \left\|g(\cdot,\cdot)\right\|_{L^{1}_v L^{p}_z} \left\|f(\cdot,\cdot) \right\|_{L^{p}_v L^{\infty}_z},
    \end{split}
    \end{equation}
    where the constant $C^+_{R,L,I_z,d,p}(B)$ depends on the collision kernel $B$, domain of random variable $I_z$, parameters $R$, $L$, $d$ and $p$.
    \begin{equation}\label{QLLp}
    \begin{split}
    \left\|Q^{R,-}(g,f)(\cdot,\cdot)\right\|_{L^p_{v,z}} \leq C^-_{R,L,I_z,d,p}(B) \|g(\cdot,\cdot)\|_{L^1_v L^{\infty}_z} \|f(\cdot,\cdot)\|_{L^p_{v,z}},\\[4pt]
    \left\|Q^{R,-}(g,f)(\cdot,\cdot)\right\|_{L^p_{v,z}} \leq C^-_{R,L,I_z,d,p}(B) \left\|g(\cdot,\cdot)\right\|_{L^{1}_v L^{p}_z} \left\|f(\cdot,\cdot) \right\|_{L^{p}_v L^{\infty}_z},
    \end{split}
    \end{equation}
    where the constant $C^-_{R,L,I_z,d,p}(B)$ depends on the collision kernel $B$, domain of random variable $I_z$, parameters $R$, $L$, $d$ and $p$.\\
    Furthermore, the complete collision operator $Q^{R}(g,f)$ satisfies
    \begin{equation}\label{QLp}
    \begin{split}
     \left\|Q^{R}(g,f)(\cdot,\cdot)\right\|_{L^p_{v,z}} \leq C_{R,L,I_z,d,p}(B) \left\|g(\cdot,\cdot)\right\|_{L^{1}_v L^{\infty}_z} \left\|f(\cdot,\cdot) \right\|_{L^p_{v,z}},\\[4pt]
     \left\|Q^{R}(g,f)(\cdot,\cdot)\right\|_{L^p_{v,z}} \leq C_{R,L,I_z,d,p}(B) \left\|g(\cdot,\cdot)\right\|_{L^{1}_v L^{p}_z} \left\|f(\cdot,\cdot) \right\|_{L^{p}_v L^{\infty}_z},
    \end{split}
    \end{equation}
    where the constant $C_{R,L,I_z,d,p}(B)$ depends on $C^+_{R,L,I_z,d,p}(B)$ and $C^-_{R,L,I_z,d,p}(B)$.\\
\end{proposition}

\begin{proof}
    The estimate of the gain term $Q^{R,+}$ replies on the definition of $L^p_{v,z}$ by duality, while the estimate of loss term $Q^{R,-}$ is achieved by capturing its convolutional structure. The complete proof can be found in Appendix~\ref{App-QLP}.
\end{proof}

\begin{remark}
  (i) As mentioned, the $L^p_{v,z}$-estimate will play an essential role in the proof of the well-posedness of either theoretical solution \eqref{ABE} or numerical solution to \eqref{PFSNK}. Note that the estimate \eqref{QLp} is currently not sufficient if we want to use the fixed point theorem to show the well-posedness (see Proposition \ref{localexistence}), due to the inconsistency of the norms on both sides of \eqref{QLp}, unless one choose to work in $L^{\infty}_{v,z}$ by selecting $p=\infty$. However, $L^{\infty}_{v,z}$ framework does not work well in proving the well-posedness of semi-discretized system \eqref{PFSNK}, as we have to rely on the Schwarz-type inequality for the projection \cite[Chapter 3]{Xiu}, i.e., $\|\mP^K_N Q^R \| \leq \|Q^R\|$, which is not always true for $L^{\infty}_{v,z}$ norm (one can take the Gibbs phenomenon as a counterexample). 
  
  (ii) To overcome the aforementioned difficulties, it suffices to present the following Corollary \ref{sumQz}. The underlying motivations include: first of all, to control the projection $\|\mP^K_N Q^R \|$ via Schwarz-type inequality, i.e., $\|\mP^K_N Q^R\| \leq \|Q^R\|$, we need to select $p=2$ in \eqref{QLp}; furthermore, to handle the existence of $L^{\infty}_z$ on the right-hand side of \eqref{QLp}, we choose to work in the higher-order Sobolev space $H^1_z$ by taking advantage of the Sobolev embedding $H^{s}_z(I_z) \hookrightarrow L^{\infty}_z(I_z), s>\frac{d_z}{2}$, which, though, restricts in the case of $z \in I_z \subset \mathbb{R}^{d_z}$ with $d_z = 1$ throughout the rest of this paper (our methodology can be generalized to the case of $d_z > 1$, where one needs to work in the higher-order space $H^s_z$ instead of $H^1_z$); besides, once considering the high-regularity in $z$, one has to apply the Leibniz rule, where the symmetric property of $Q^R$ will be helpful to deal with higher-order derivative of function, which is why we show two inequalities in \eqref{QLp}.
\end{remark}

\begin{corollary}\label{sumQz}
    Let the collision kernel $ B $ satisfy the assumptions \eqref{kernel}-\eqref{cutoff}, \eqref{kinetic}-\eqref{Assump_B} and truncation parameters $R$ and $L$ satisfy \eqref{RL}, then, for $z \in I_z \subset \mathbb{R}^{d_z}$ with $d_z=1$, the truncated collision operators $ Q^{R}(g,f)$ satisfies the following estimates:
    \begin{equation}\label{QL2H1}
        \| Q^R(g,f)(\cdot,\cdot) \|_{L^2_v H^1_z} \leq C_{0}(B) \left\|g(\cdot,\cdot)\right\|_{L^{1}_vH^1_z} \left\|f(\cdot,\cdot) \right\|_{L^{2}_vH^1_z},
    \end{equation}
    \begin{equation}\label{Q_HkZ}
        \| Q^R(g,f)(\cdot,\cdot) \|_{H^k_v H^1_z} \leq C_{k}(B) \left\|g(\cdot,\cdot)\right\|_{H^k_v H^1_z} \left\|f(\cdot,\cdot) \right\|_{H^k_v H^1_z},
    \end{equation}
    where $C_{0}(B)$ and $C_{k}(B)$ are constants depending on the collision kernel $B$, domain of random variable $I_z$, parameters $R$, $L$, $d$ and $p$.
\end{corollary}

\begin{proof}
By selecting $p=2$ in \eqref{QLp}, it yields that
\begin{equation}\label{Qgf1}
    \left\|Q^{R}(g,f)\right\|_{L^2_{v,z}} \leq C_{R,L,I_z,d,2}(B) \left\|g\right\|_{L^{1}_v L^{\infty}_z} \left\|f \right\|_{L^2_{v,z}},
\end{equation}
\begin{equation}\label{Qgf2}
    \left\|Q^{R}(g,f)\right\|_{L^2_{v,z}} \leq C_{R,L,I_z,d,2}(B) \left\|g\right\|_{L^{1}_v L^{2}_z} \left\|f\right\|_{L^{2}_v L^{\infty}_z}.
\end{equation}

To prove the estimate \eqref{QL2H1}, we apply $\partial_z$ to $Q^R(g,f)$ and consider the Leibniz rule, 
\begin{equation}\label{Q_z}
\partial_z Q^R(g,f) = Q^R(\partial_z g, f) + Q^R(g ,\partial_z f) + Q_{B^{*}}^R(g,f), 
\end{equation}
where $Q_{B^{*}}^R$ is defined by substituting $\partial_z B_T$ by $B_T$ in $Q^R$.
Therefore, 
\begin{equation}\label{QL2H1proof}
    \begin{split}
        \| Q^R(g,f)\|_{L^2_v H^1_z}^2 =&  \|Q^R(g,f)\|_{L^2_{v,z}}^2 + \|\partial_z Q^R(g,f)\|_{L^2_{v,z}}^2 \\[4pt]
        \leq & \|Q^R(g,f)\|_{L^2_{v,z}}^2 + \|Q^R(\partial_z g, f)\|_{L^2_{v,z}}^2 + \|Q^R(g,\partial_z f)\|_{L^2_{v,z}}^2 +  \|Q_{B^{*}}^R(g,f)\|_{L^2_{v,z}}^2 \\[4pt]
        \leq & C^2_{R,L,I_z,d,2}(B) \left( \left\|g\right\|_{L^{1}_v L^{\infty}_z}^2 \left\|f \right\|_{L^2_{v,z}}^2 + \left\|\partial_z g\right\|_{L^{1}_v L^{2}_z}^2 \left\|f\right\|_{L^{2}_v L^{\infty}_z}^2 + \left\|g\right\|_{L^{1}_v L^{\infty}_z}^2 \left\|\partial_z f \right\|_{L^2_{v,z}}^2 \right)\\[4pt]
        & + C^2_{R,L,I_z,d,2}(B^*) \left\|g\right\|_{L^{1}_v L^{\infty}_z}^2 \left\|f \right\|_{L^2_{v,z}}^2 \\[4pt]
        \leq & C^2_{R,L,I_z,d,2}(B) \left\|g\right\|_{L^{1}_vH^1_z}^2 \left\|f \right\|_{L^{2}_vH^1_z}^2 + C^2_{R,L,I_z,d,2}(B^*) \left\|g\right\|_{L^{1}_vH^1_z}^2 \left\|f \right\|_{L^{2}_vH^1_z}^2 \\[4pt]
        \leq & C^2_{0}(B) \left\|g\right\|_{L^{1}_vH^1_z}^2 \left\|f\right\|_{L^{2}_vH^1_z}^2,
    \end{split}
\end{equation}
where $C^2_{0}(B)$ is chosen as $\max\{C^2_{R,L,I_z,d,2}(B), C^2_{R,L,I_z,d,2}(B^*) \}$, and in the second inequality we apply \eqref{Qgf1} to estimate $\|Q^R(\partial_z g, f)\|_{L^2_{v,z}}^2$, while \eqref{Qgf2} is used for $\|Q^R(g, \partial_z 
 f)\|_{L^2_{v,z}}^2$. In the third inequality, the Sobolev embedding $H^{s}_z(I_z) \hookrightarrow L^{\infty}_z(I_z)$ is utilized for $s > \frac{d_z}{2}$ with $z \in I_z \subset \mathbb{R}^{d_z}$ and $d_z = 1$.
Hence, \eqref{QL2H1} is proved by taking the square root of \eqref{QL2H1proof} above.

Furthermore, to obtain the estimate \eqref{Q_HkZ}, one needs to apply \eqref{QL2H1} as well as the Leibniz rule in $v$ when considering the high-order derivatives in velocity. The detailed proof can directly follow the deterministic case \cite[Proposition 3.2]{HQY21}. 
\end{proof}

\subsubsection{Mass conservative property with uncertainties}
\label{subsubsec:mass}

The Fourier spectral method for velocity discretization, albeit efficient, exhibits a limitation in its preservation of some crucial physical properties, such as positivity. Nonetheless, a pivotal attribute, namely the conservation of mass, is preserved throughout time, thereby offering valuable control over the numerical solution $f_N^K$. To articulate this point precisely, we present the following Lemma by following the deterministic counterpart in \cite[Lemma 2.1]{HQY21}:

\begin{lemma}\label{lemma:conv}
	The numerical system \eqref{PFSNK} preserves mass in the sense that, for $r \geq 0$,
	\begin{equation} \label{mass_conserve}
	\| \int_{\D} f_N^K(t,v,\cdot)  \,\rd v \|_{H_z^r} = \| \int_{\D} f^{0,K}(v,\cdot)  \,\rd v \|_{H_z^r}.
	\end{equation}
\end{lemma}
\begin{proof}
    Consider that
    \begin{equation}\label{m1}
    \begin{split}
        \| \int_{\D} f_N^K(t,v,\cdot)  \,\rd v \|_{H_z^r}  = &\sum_{|\mathbf{n}|=0}^N \sum_{|\mathbf{k}|=0}^K  \| f_{\mathbf{n}}^{\mathbf{k}}(t) \Psi^{\mathbf{k}}(\cdot) \|_{H_z^r} \int_{\D} \e^{\im \frac{\pi}{L} \mathbf{n}\cdot v} \,\rd v \\
        =& (2L)^d \sum_{|\mathbf{k}|=0}^K  \| f_{\mathbf{0}}^{\mathbf{k}}(t) \Psi^{\mathbf{k}}(\cdot) \|_{H_z^r},
    \end{split}
    \end{equation}
    where $f_{\mathbf{0}}^{\mathbf{k}}(t)$ is the $0$-th mode (in velocity) of the numerical solution $f_N^K$. Furthermore, $f_{\mathbf{0}}^{\mathbf{k}}(t)$ satisfies the equation
    \begin{equation}
        \frac{\rd }{\rd t} f_{\mathbf{0}}^{\mathbf{k}}(t) = Q^{R,\kk}_{\mathbf{0}}(t),
    \end{equation}
    in fact, $Q^{R,\kk}_{\mathbf{0}}(t) \equiv 0 $, one can see this from the definition of $Q^{R,\kk}_{\n}$ in \eqref{QRnk}, since $G(\mathbf{l}, \mathbf{m}) \equiv 0 $ when $\n = \mathbf{l} + \mathbf{m} = \mathbf{0}$. 

    This implies that $f_{\mathbf{0}}^{\mathbf{k}}(t)$ remains to be constant as time evolves, i.e., $f_{\mathbf{0}}^{\mathbf{k}}(t) = f_{\mathbf{0}}^{0,\mathbf{k}}$ such that, following \eqref{m1},
    \begin{equation}
        \| \int_{\D}  f_N^K(t,v,\cdot)  \,\rd v \|_{H_z^r} = (2L)^d \sum_{|\mathbf{k}|=0}^K  \| f_{\mathbf{0}}^{0,\mathbf{k}} \Psi^{\mathbf{k}}(\cdot) \|_{H_z^r} = \| \int_{\D} f^{0,K}(t,v,\cdot)   \,\rd v \|_{H_z^r}.
    \end{equation}
    Hence, \eqref{mass_conserve} is proved.
    Note that, for the initial approximation $f^{0,K}$, it is typically obtained by applying the projection $\mathcal{P}^K$ on $f^0$, therefore, we particularly have, for $r=1$,
    \begin{equation}\label{mass_conse_PK}
    \begin{split}
         \| \int_{\D}  f_N^K(t,v,\cdot)  \,\rd v \|_{H^1_z}  = \| \int_{\D} f^{0,K}(v,\cdot)   \,\rd v \|_{H^1_z} =& \| \int_{\D} \mathcal{P}^K f^{0}(v,\cdot)   \,\rd v \|_{H^1_z}\\
         =& \| \int_{\D}  f^{0}(v,\cdot) \,\rd v \|_{H^1_z},
    \end{split}
    \end{equation}
    which is an initial quantity $\| f^{0}(\cdot,\cdot) \|_{L^1_v H^1_z} $ denoted by $E_1^{f^0}$ through the rest of the paper.
\end{proof}


\section{Well-posedness of the numerical solution \texorpdfstring{$f_N^K$}{fN} with uncertainties}
\label{sec:main}

In this section, we establish the well-posedness of the numerical solution $f_N^K$ solved from the semi-discretized system \eqref{PFSNK} on any arbitrarily given time interval $[0, T]$. It is noteworthy that, akin to the deterministic scenario discussed in \cite{HQY21}, the primary challenge in the proof stems from the fact that the numerical solution $f_N^K$ may not be inherently non-negative. This characteristic arises due to the intrinsic limitation of the spectral projection $\mP_N^K$, even though the analytic counterpart $f$ to the original problem \eqref{IBE} is always non-negative.

Similar to our preceding work \cite{LQ2022}, our strategy here is to extend the local well-posedness result until covering up the whole prescribed time interval $[0, T]$, by applying the delicate energy estimate. In Section \ref{subsec:propagation}, we obtain the propagation of $H_v^k H^r_z$ norm of the solution $f_N^K$ with uncertainties, where the key is a ``good" control of the negative part of $f_N^K$, aside from using the mass conservation, i.e., Lemma \ref{lemma:conv} as well as {\it a priori} $L^1_v H^1_z$ bound. In Section \ref{subsec:wellposed}, the local existence and uniqueness can be shown via a fixed point theorem as long as the time interval $[0, \tau]$ is chosen to be sufficiently small, where the negative part of the numerical solution $f_N^K$ is under control within the same period of time with large enough $N$ and $K$, which in turn implies that the initial $L^1_v H^1_z$ bound of $f_N^K$ can be preserved at time $\tau$. Therefore, the same procedure can be repeated iteratively to extend the solution up to the prescribed final time $T$, where the values of designated parameters $N$, $K$ and $\tau$ are shown to be the same in each iteration as the beginning one.

\subsection{Propagation of Sobolev estimates of numerical solution with uncertainties}
\label{subsec:propagation}

\subsubsection{\texorpdfstring{$H_v^kH^r_z$}{Hvkrz} estimate of \texorpdfstring{$f_N^K$}{fNK}}
\label{subsub:Hk}

We now consider the propagation of the $H_v^k H_z^r$-norm of $f_N^K$.


\begin{proposition}\label{regularity}
Let the truncation parameters $R$, $L$ satisfy \eqref{RL} and assume the collision kernel $B$ satisfies \eqref{kernel}--\eqref{cutoff} and \eqref{kinetic}--\eqref{Assump_B}.
For the semi-discretized system \eqref{PFSNK}, assume that the initial condition $f_{N}^{0,K}(v,z) \in H^k_v H_z^r $ for some integers $k \geq 0$ and $r \geq 1$, i.e., 
\begin{equation*}
    \| f_{N}^{0,K}(\cdot,\cdot)\|_{H^k_v H_z^r} \leq D^{f^0}_{k,r}.
\end{equation*}
If $f_N^K(t,v,z)$ is further assumed to have a $L_v^1 H_z^1$ bound up to some time $t_0$, i.e., 
\begin{equation}\label{Priori_fZ}
\forall t\in [0,t_0], \quad \left\| f_N^K(t,\cdot,\cdot) \right\|_{L_v^1 H_z^1} \leq  2E_1^{f^0}, 
\end{equation}
then there exists a $K_{k,r}(t_0) > 0$ depending on $t_0$, $k$, $r$, $E_1^{f^0}$, $D^{f^0}_{k,r}$ and $C_{0}(B)$ such that,
\begin{equation}\label{fN_HkZ}
\forall t\in[0,t_0], \quad \left\| f_N^K(t,\cdot,\cdot) \right\|_{H^k_v H^r_z} \leq K_{k,r}(t_0).
\end{equation}
\end{proposition}


\begin{proof} 
We start with defining 
\begin{equation*}
\begin{split}
    Q_{B^{l}}^R(h,f)(v,z) =& \int_{\mathbb R^d}\int_{\mathbb{S}^{d-1}} \partial_z^l B_T(\sigma\cdot \hat{q},|q|,z) \left[h(v_*',z)f(v',z) - h(v_{*},z)f(v,z)\right] \,\rd{\sigma}\, \rd{v_{*}}\\[4pt]
    :=& \int_{\mathbb R^d}\int_{\mathbb{S}^{d-1}} \partial_z^l B_T \left[h'_{*} f' - h f \right]\,\rd{\sigma}\, \rd{v_{*}}\,.
\end{split}
\end{equation*}
where, considering \eqref{Assump_B}, we obtain that all $l$-th order $z$-derivative of $B$ is uniformly bounded, i.e., $| \partial_z^l B_T | \leq C_B$ for all $z \in I_z$ and $|l|\leq r$. 
Recalling \eqref{QR}, the only difference between $Q_{B^{l}}^R$ and $Q^R$ is that $B_T$ is substituted by $\partial_z^l B_T$. 
Using the Leibniz rule to the collision operator in $v$ and $z$ yields that
\begin{equation}
\begin{split}
\quad\partial_v^{\nu}\partial_z^l Q^R(h, f) & =
\sum_{|\mu|=0}^{|\nu|}\binom{\nu}{\mu} \sum_{|n|=0}^{|l|-1}\sum_{|m|=0}^{|n|} \binom{l}{n}\binom{n}{m} Q_{B^{l-n}}^R (\partial_v^{\mu}\partial_z^m h, \partial_v^{\nu-\mu} \partial_z^{n-m}f) \\[4pt]
& \quad + \sum_{|\mu|=0}^{|\nu|}\binom{\nu}{\mu} \sum_{|m|=1}^{|l|-1}\binom{l}{m}Q^R(\partial_v^{\mu}\partial_z^m h, \partial_v^{\nu-\mu}\partial_z^{l-m}f)  \\[4pt]
& \quad + \sum_{|\mu|=0}^{|\nu|}\binom{\nu}{\mu} Q^R(\partial_v^{\mu}h, \partial_v^{\nu-\mu}\partial_z^l f) 
+ \sum_{|\mu|=0}^{|\nu|}\binom{\nu}{\mu} Q^R(\partial_v^{\mu}\partial_z^l h, \partial_v^{\nu-\mu}f). 
\end{split}
\end{equation}
Substituting $h$ and $f$ by $f_N^K$, we derive the following estimate: 
\begin{equation}\label{Q_Leibniz2}
\begin{split}
& \| \partial_z^l Q^R(\fe^K, \fe^K) \|_{H_v^k L_z^2} =
  \sum_{|\nu|=0}^{k}\|\partial_v^{\nu}\partial_z^l Q^R(\fe^K, \fe^K)\|_{L^2_{v,z}} \\[4pt]
\leq& C_{R,L,I_z,d,2}(B) \sum_{|\nu|=0}^{k}\sum_{|\mu|=0}^{|\nu|}\binom{\nu}{\mu} \sum_{|n|=0}^{|l|-1}\sum_{|m|=0}^{|n|} \binom{l}{n}\binom{n}{m} 
\|\partial_v^{\mu}\partial_z^m \fe^K\|_{L^2_{v} H^1_z} \,
\|\partial_v^{\nu-\mu}\partial_z^{n-m}\fe^K \|_{L^2_{v,z}}  \\[4pt]
& + C_{R,L,I_z,d,2}(B) \sum_{|\nu|=0}^{k}\sum_{|\mu|=0}^{|\nu|}\binom{\nu}{\mu} \sum_{|m|=1}^{|l|-1}\binom{l}{m} \|\partial_v^{\nu}\partial_z^m \fe^K \|_{L^2_{v}H^1_z} \, \|\partial_v^{\nu-\mu}\partial_z^{l-m} \fe^K\|_{L^2_{v,z}} \\[4pt]
& + C_{R,L,I_z,d,2}(B) \sum_{|\nu|=0}^{k}\sum_{|\mu|=0}^{|\nu|}\binom{\nu}{\mu} 
\|\partial_v^{\mu} \fe^K \|_{L^2_v L_z^{\infty}} \, \|\partial_v^{\nu-\mu}\partial_z^l \fe^K \|_{L^2_{v,z}} \\[4pt]
& + C_{R,L,I_z,d,2}(B) \sum_{|\nu|=0}^{k}\sum_{|\mu|=0}^{|\nu|}\binom{\nu}{\mu} \|\partial_v^{\mu}\partial_z^l \fe^K \|_{L^2_{v,z}} \, \|\partial_v^{\nu-\mu} \fe^K \|_{L^2_v L_z^{\infty}}, 
\end{split}
\end{equation}
where \eqref{Qgf1} is used in the inequality above and we further apply the embedding relation. 

Then, we claim the following statement holds,
\begin{equation}\label{PZ_fNK} 
\forall t\in [0,t_0], \qquad \| \partial_z^l f_N^K(t,\cdot,\cdot)\|_{H_v^k L_z^2} \leq K_{k,l}(t_0), \quad  0 \leq |l| \leq r,
\end{equation}
and will use mathematical induction to show it.
When $r=0$, based on the assumption \eqref{Priori_fZ}, we can go through almost the same analysis as in the deterministic problem \cite{HQY21} (except that we need to take $H_z^1$ norm on both sides of the equation \eqref{PFSNK} and employ \eqref{QL2H1}), then 
\begin{equation}\label{fNK_Det} 
\forall t\in [0,t_0], \qquad \|f_N^K(t,\cdot,\cdot)\|_{H_v^k H_z^1} \leq K_{k,1}(t_0). 
\end{equation}
This further implies that 
    \begin{equation}\label{fNK_Inf}
    \|f_N^K(t,\cdot,\cdot) \|_{H_v^k L_z^2} \leq C_{I_z} \|f_N^K(t,\cdot,\cdot) \|_{H_v^k L_z^{\infty}} \leq \|f_N^K(t,\cdot,\cdot) \|_{H_v^k H^1_z} \leq K_{k,1}(t_0).
    \end{equation}
Therefore, \eqref{PZ_fNK} holds for $r=0$. 

Now we are in position to show that, when $|l| = r$, $\left\| \partial_z^l \fe^K(t,\cdot,\cdot) \right\|_{H_v^k L^2_z} \leq K_{k,r}(t_0)$ by presuming that $\left\| \partial_z^l \fe^K(t,\cdot,\cdot) \right\|_{H_v^k L^2_z} \leq K_{k,r-1}(t_0)$ for $0 \leq |l| \leq r-1$.

Taking $\partial_v^{\nu}\partial_z^l$ on both sides of \eqref{PFSNK}, 
\begin{equation}\label{Derivative_vz}
	\partial_t \partial_v^{\nu}\partial_z^l \fe^K(t,v,z) = \mathcal{P}_N^K \partial_v^{\nu} \partial_z^l Q^R(\fe^K, \fe^K)(t,v,z),
\end{equation}
multiplying both sides of \eqref{Derivative_vz} by $\partial_v^{\nu}\partial_z^l \fe^K$, and integrating over $v \in \mathcal{D}_L $ and $ z \in I_z$ gives 
\begin{equation}\label{fNK2}
\begin{split}
\frac{1}{2} \frac{\rd}{\rd t} \| \partial_v^{\nu} \partial_z^l \fe^K\|_{L_{v,z}^2}^2  =& \int_{I_z} \int_{\D} \mathcal{P}_N^K  \partial_v^{\nu}\partial_z^l Q^R(\fe^K, \fe^K) \,
\partial_v^{\nu}\partial_z^l \fe^K \,\rd v \,\pi(z) \,\rd z \\
\leq & \|\mathcal{P}_N^K \partial_v^{\nu}\partial_z^l Q^R(\fe^K, \fe^K) \|_{L_{v,z}^2}\, 
\|\partial_v^{\nu}\partial_z^l \fe^K \|_{L_{v,z}^2 } \\[4pt]
\leq & \| \partial_v^{\nu}\partial_z^l Q^R(\fe^K, \fe^K) \|_{L_{v,z}^2}\, \|\partial_v^{\nu}\partial_z^l \fe^K \|_{L_{v,z}^2}.  
\end{split}
\end{equation}
Adding up \eqref{fNK2} for all $|\nu|\leq k$ and using the Cauchy-Schwarz inequality, one derives that
\begin{equation}\label{fNK_Inq}
\frac{1}{2} \frac{\rd}{\rd t} \| \partial_z^l \fe^K \|_{H_v^{k}L_z^2}^2
\leq \| \partial_z^l Q^R(\fe^K, \fe^K)\|_{H_v^{k}L_z^2} \|\partial_z^l \fe^K \|_{H_v^{k}L_z^2}. 
\end{equation}
Therefore, for $|l|=r$,
\begin{equation}\label{fe_K}
\begin{split}
 \frac{1}{2} \frac{\rd}{\rd t} \| \partial_z^l \fe^K \|_{H_v^{k}L_z^2} & \leq
  \| \partial_z^l Q^R(\fe^K, \fe^K)\|_{H_v^{k}L_z^2}  \\[4pt]
& = \| \partial_z^l Q^R(\fe^K, \fe^K)\|_{H_v^{k-1}L_z^2} + \sum_{|\nu|=k}
\| \partial_v^{\nu}\partial_z^l Q^R(\fe^K, \fe^K) \|_{L^2_{v,z}}. 
\end{split}
\end{equation} 
For the first term on the right-hand side of \eqref{fe_K}, from \eqref{Q_Leibniz2} one gets that
\begin{equation}\label{dd1}
\begin{split}
& \| \partial_z^l Q^R(\fe^K, \fe^K)\|_{H_v^{k-1}L_z^2} \\[4pt] 
\leq & \| \partial_z^l Q^R(\fe^K, \fe^K)\|_{H_v^{k}L_z^2} \\[4pt]
\leq & C_{R,L,I_z,d,2}(B)  \Big\{ C_{r-1,k} \|\fe^K\|_{H_v^k H_z^{r-1}}^2 + 
 C_k \|\partial_z^l \fe^K \|_{H_v^k L_z^2} \|\fe^K\|_{H_v^k L_z^\infty} \\[4pt]
 & \qquad \qquad \qquad \qquad + 
 C_k^{'}  \|\partial_z^l \fe^K \|_{H_v^k L^2_z}
 \left(  \| \fe^K \|_{H_v^k L_z^2} + \| \partial_z \fe^K \|_{H_v^k L^2_z} \right)  \Big\} \\[4pt]
\leq & C_{k,\tilde K,r-1}(B) + C_k(B)
 \| \partial_z^l \fe^K \|_{H_v^k L^2_z}  \left( \| \fe^K \|_{H_v^k L_z^{\infty}} +  \| \fe^K \|_{H_v^k L_z^2} + \| \partial_z \fe^K \|_{H_v^k L^2_z}  \right) \\[4pt]
\leq & C_{k,\tilde K,r-1}(B) + C_{k,K_{k,0}}(B)
 \| \partial_z^l \fe^K \|_{H_v^k L^2_z}, 
\end{split}
\end{equation}
where $C_{k,K_{k,0}(B)} > 0$ depends on $C_{R,L,I_z,d,2}(B)$, $k$ and $K_{k,0}$; while $C_{k,\tilde K,r-1}(B)> 0$ depends on $C_{R,L,I_z,d,2}(B)$, $k$, $r-1$ and $\tilde K:=\max_{|l|\leq r-1} K_{k,|l|}$. The last two terms in the second inequality arise from taking $|m|=|n|=|l|-1$ and $m=l$ in the first and second terms in \eqref{Q_Leibniz2} respectively. The induction hypothesis, embedding relation $ H^1_z \hookrightarrow L^{\infty}_z$, and \eqref{fNK_Inf} are used in the last inequality.\\
For the second term on the right-hand side of \eqref{fe_K}, similar to \eqref{Q_Leibniz2} except that $|\nu|=k$ now, we can split the summation into three cases: $|\mu| = 0$, $0 < |\mu| \leq k$, and $|\mu|=k$, then it yields that
\begin{small}
\begin{equation}\label{dd2}
\begin{split}
& \sum_{|\nu|=k}\| \partial_v^{\nu}\partial_z^l Q^R(\fe^K, \fe^K) \|_{L^2_{v,z}}\\[4pt]  
\leq & C_{R,L,I_z,d,2}(B) \Big\{ \widetilde C_{r-1,k}\|\fe^K\|_{H_v^{k-1}H_z^{r-1}} + 2 \sum_{|\nu|=k} \left( \|\partial_v^{\nu}\fe^K\|_{L^2_v L_z^{\infty}} \|\partial_z^l \fe^K \|_{L^2_{v,z}}
+ \| \fe^K\|_{L^2_v L_z^{\infty}} \|\partial_v^{\nu}\partial_z^l\fe^K\|_{L^2_{v,z}} \right) \\
& \qquad\qquad\qquad \qquad  +  C_k^{'}  \|\partial_z^l \fe^K \|_{H_v^k L^2_z}
 \left(  \| \fe^K \|_{H_v^k L_z^2} + \| \partial_z \fe^K \|_{H_v^k L^2_z} \right) \Big\} \\[4pt]
\leq & \widetilde C_{k,\tilde K,r-1}(B) + 2 C_{R,L,I_z,d,2}(B) \|\partial_z^l \fe^K\|_{H_v^k L_z^2}\left(
 \| \fe^K \|_{H_v^k L_z^{\infty}} + \| \fe^K \|_{H_v^k L_z^2} + \| \partial_z \fe^K \|_{H_v^k L^2_z} \right)  \\[4pt]
 \leq & \widetilde C_{k,\tilde K,r-1}(B) + \widetilde C_{k,K_{k,0}}(B) \|\partial_z^l \fe^K\|_{H_v^k L_z^2}, 
\end{split}
\end{equation}
\end{small}
where $\widetilde C_{k,\tilde K,l-1}(B) > 0$ is a constant depending on $C_{R,L,d,2}(B)$, $k$, $l-1$ and $\widetilde K$; while $\widetilde C_{k,K_{k,0}}(B)$ is a constant depending on $C_{R,L,d,2}(B)$, $k$ and $K_{k,0}$. The induction hypothesis, embedding relation $ H^1_z \hookrightarrow L^{\infty}_z$ and \eqref{fNK_Inf} are used in the last inequality above.

By substituting the estimates \eqref{dd1} -\eqref{dd2} into \eqref{fe_K} and applying the Gr\"onwall's inequality, we find, for $t \in [0, t_0]$,
\begin{equation}
\begin{split}
\|\partial_z^l \fe^K \|_{H_v^k L^2_z} \leq & \e^{t \bar C_{k,K_{k,0}}(B)} \left( \|\partial_z^l \fe^{0,K} \|_{H_v^k L^2_z} + \bar C_{k,\tilde K,r-1}(B) \right) \\[4pt]
\leq & \e^{t \bar C_{k,K_{k,0}}(B)}
\left( D_{k,r}^{f^0} + \bar C_{k,\tilde K,r-1}(B) \right) := \bar{K}_{k,r}(t_0), 
\end{split}
\end{equation}
where the constant $C_{k,K_{k,0}}(B) > 0$ depends on $C_{k,K_{k,0}}(B)$ and $\widetilde C_{k,K_{k,0}}(B)$, and the constant $\bar C_{k,\tilde K,r-1}(B) > 0$ depends on $ C_{k,\tilde K,l-1}(B)$ and $\widetilde C_{k,\tilde K,l-1}(B)$.
Thus, by mathematical induction, the statement \eqref{PZ_fNK} is proved. 

Finally, by adding up all $l$ for $0\leq |l|\leq r$, we can obtain that 
$$ \forall t\in [0,t_0], \qquad \|\fe^K(t,\cdot,\cdot)\|_{H_v^k H_z^r}\leq K_{k,r}(t_0), $$
where $K_{k,r}(t_0) := \sum_{|l|=0}^r \bar{K}_{k,|l|}(t_0) > 0$ depends on $k$, $r$, $E_1^{f^0}$, $D^{f^0}_{k,r}$ and $C_{0}(B)$, and will be increased as $t_0$ rises.
\end{proof}


\subsubsection{\texorpdfstring{$L_v^2 H^1_z$}{Lv2H1z} estimate of the negative part \texorpdfstring{$f^{K,-}_N$}{f-NK}}
\label{subsub:L2H1f-}

We now proceed to estimate the negative part of $f_N^K$, i.e., $f^{K,-}_N$, which will play a key role in extending the local well-posedness of $f_N^K$ to a long prescribed time interval.
The result can be summarized into the following Proposition: 

\begin{proposition}\label{Pro_fnk_neg}
Let the truncation parameters $R$, $L$ satisfy \eqref{RL} and assume the collision kernel $B$ satisfies \eqref{kernel}--\eqref{cutoff} and \eqref{kinetic}--\eqref{Assump_B}. 
For the semi-discretized system \eqref{PFSNK}, assume that the initial condition $f_{N}^{0,K}(v,z) \in H^1_v H_z^1 $, i.e., 
\begin{equation*}
    \| f_{N}^{0,K}(\cdot,\cdot)\|_{H^1_v H_z^1} \leq D^{f^0}_{1,1},
\end{equation*}
and $f_N^K(t,v,z)$ has a $L_v^1 H_z^1$ bound up to some time $t_0$, i.e., 
\begin{equation}\label{Priori_fZ_2}
\forall t\in [0,t_0], \quad \left\| f_N^K(t,\cdot,\cdot) \right\|_{L_v^1 H_z^1} \leq  2E_1^{f^0}. 
\end{equation}
Then, the following estimate holds for the numerical solution $f_N^K$,
\begin{equation}\label{FNK0111}
    \forall t\in [0,t_0], \quad \|f_N^K(t,\cdot,\cdot)\|_{L_v^2 H_z^1} \leq K_{0,1}(t_0), \quad \|f_N^K(t,\cdot,\cdot)\|_{H_v^1 H_z^1} \leq K_{1,1}(t_0),
\end{equation}
and for the negative part of the numerical solution $f_N^{K,-}$,
\begin{equation}\label{fe_minus}
\forall t \in [0,t_0], \quad  \left\| f_N^{K,-} (t,\cdot,\cdot) \right\|_{L^2_v H^1_z}
 \leq \e^{ t \, \mathcal{C}_{K_{0,1}(t_0)}} \left( \| f_N^{0,K,-}(\cdot, \cdot) \|_{L^2_v H^1_z} + \frac{\mathcal{C}_{K_{1,1}(t_0)}}{N} \right), 
\end{equation}
where $\mathcal{C}_{K_{0,1}(t_0)} > 0$ depends on $C_0(B)$, $E_1^{f^0}$ and $K_{0,1}(t_0)$, while $\mathcal{C}_{K_{1,1}(t_0)}>0$ depends on $C_0(B)$ and $K_{1,1}(t_0)$.
\end{proposition}

\begin{proof}
Firstly, \eqref{FNK0111} can be directly derived from Proposition \ref{regularity} by selecting $k=0, r=1$ and $k=1,r=1$, respectively.

To estimate $f_N^{K,-}$ in $\|\cdot\|_{L^2_v H^1_z}$, we take $\partial_z^l$ on both sides of the numerical system \eqref{PFSNK} with $|l|=0,1$,
\begin{equation}
\partial_t \partial_z^l f_N^K(t,v,z) = \partial_z^l \mathcal{P}_N^K Q^R(f_N^K, f_N^K)(t,v,z),
\end{equation}
which can be further re-written as 
\begin{equation} \label{eqn-g1} 
\partial_t \partial_z^l f_N^K(t,v,z) = \partial_z^l Q^R(f_N^K, f_N^K)(t,v,z) + \partial_z^l E^{K}_N(t,v,z), 
\end{equation}
where 
\begin{equation} \label{eqn-E}
E_N^{K}(t,v,z):= \mP_N^K  Q^R(f_N^K, f_N^K)(t,v,z)  -  Q^R(f_N^K, f_N^K)(t,v,z). 
\end{equation}

Define the indicator function $\mathbf{1}_{\left\lbrace \partial_z f_N^K \leq 0\right\rbrace}$ as follows: 
\begin{equation}\label{chi1}
	\mathbf{1}_{\left\lbrace  f_N^K \leq 0 \right\rbrace} (v) := \begin{cases}
		1,& \text{for}\ v\ \text{s.t.}\  f_N^K(\cdot,v,\cdot) \leq 0, \\
		0,& \text{for}\ v\ \text{s.t.}\  f_N^K(\cdot,v,\cdot) > 0.
	\end{cases}
\end{equation}
Note that $\partial_z^l f_N^K(t,v,z) = \partial_z^l f_N^{K,+}(t,v,z) - \partial_z^l f_N^{K,-}(t,v,z)$.
After multiplying both sides of \eqref{eqn-g1} by $\mathbf{1}_{\left\lbrace \partial_z^l f_N^K \leq 0\right\rbrace}$, the left-hand side becomes 
\begin{equation}\label{gn_char} 
\begin{split}
     &\partial_z^l f_N^K  \mathbf{1}_{\left\lbrace \partial_z^l f_N^K \leq 0\right\rbrace} (t,v,z) \partial_t  \partial_z^l f_N^K(t,v,z) \\[4pt]
     =& -  \partial_z^l f_N^{K,-}(t,v,z) \partial_t \left[ \partial_z^l f_N^{K,+}(t,v,z) -  \partial_z^l f_N^{K,-}(t,v,z)\right]\\[4pt]
     =& \  \partial_z^l f_N^{K,-}(t,v,z) \partial_t  \partial_z^l f_N^{K,-}(t,v,z),
\end{split}
\end{equation}

By considering \eqref{FNK0111}, we have 
$\displaystyle \| f_N^K(t,\cdot,\cdot) \|_{H_{v}^1 H^1_z} \leq K_{1,1}$. 
Therefore, for $|l|=0,1$, 
\begin{equation}\label{EN}
\begin{split}
&\int_{I_z} \int_{\D} \partial_z^l E^K_N(t,v,z) \partial_z^l f_N^K \mathbf{1}_{\left\lbrace \, \partial_z^l f_N^K \leq 0 \right\rbrace}(t,v,z) \,\rd v \,\pi(z) \,\rd z \\[4pt]
= & - \int_{I_z} \int_{\D} \partial_z^l E^K_N(t,v,z) \, \partial_z^l f_N^{K,-}(t,v,z) \,\rd v \,\pi(z) \,\rd z\\[4pt]
\leq & \ \| E^{K}_N \|_{L_{v}^2 H_z^1}\, \| \partial_z^l f_N^{K,-}\|_{L_{v}^2 L_z^2} \\[4pt]
\leq & \ \frac{C_{\mP}C_{1}(B)}{N} \| f_N^K \|^2_{H_{v}^1 H^1_z} \, \|f_N^{K,-}\|_{L_{v}^2 H_z^1} \\[4pt]
\leq & \ \frac{C_{\mP}C_{1}(B) \, K^2_{1,1}}{N} \|f_N^{K,-}\|_{L_{v}^2 H_z^1}\,.
\end{split}
\end{equation}  
where, in the second inequality, we apply the estimate of remainder term $E^{K}_N$, 
\begin{equation}
\begin{split}
\| E^{K}_N \|_{L_v^2 H_z^1} & = \left\| \mP_N^K  Q^R(f_N^K, f_N^K)  -  Q^R(f_N^K, f_N^K) \right\|_{L_{v}^2 H_z^1} \\[4pt]
& \leq \frac{C_{\mP}}{N} \left\| Q^R(f_N^K, f_N^K) \right\|_{H_{v}^1 H^1_z}  \leq \frac{C_{\mP}C_{1}(B)}{N} \| f_N^K \|^2_{H_{v}^1 H_z^1}.
\end{split}
\end{equation}
Here we utilize the property of the projection operator $\mP_N^K$ and the estimate \eqref{Q_HkZ} in the last inequality. 

Then, we can estimate another term on the right-hand side of \eqref{eqn-g1}, for $|l|=0,1$,
\begin{equation}\label{Q_g}
\begin{aligned}
& \int_{I_z} \int_{\D} \partial_z^l Q^R(f_N^K, f_N^K)(t,v,z)  \, \partial_z^l f_N^K \mathbf{1}_{\left\lbrace \partial_z^l f_N^K \leq 0 \right\rbrace}(t,v,z) \,\rd v \,\pi(z) \,\rd z \\[4pt]
= & - \int_{I_z} \int_{\D} \partial_z^l Q^R(f_N^K, f_N^K)(t,v,z) \, \partial_z^l f_N^{K,-}(t,v,z) \,\rd v \,\pi(z) \,\rd z \\[4pt]
\leq &\  \| Q^R(f_N^K, f_N^K)\|_{L_{v}^2 H^1_z} \| f_N^{K,-}\|_{L_{v}^2H^1_z}. 
\end{aligned}
\end{equation}  

By the Leibniz rule as in \eqref{Q_Leibniz2}, note that, for $|l|=1$,
\begin{equation}\label{QR_Z_split}
 \partial_z^l Q^R(f_N^K, f_N^K)(t,v,z) = Q^R(\partial_z^l f_N^K, f_N^K)(t,v,z) + Q^R(f_N^K, \partial_z^l f_N^K)(t,v,z), 
\end{equation}    
Let us now look at the term on the right-hand side of \eqref{QR_Z_split}, namely
\begin{equation*}
    Q^R(\partial_z^l f_N^K, f_N^K) = Q^{R,+}(\partial_z^l f_N^K, f_N^K) - Q^{R,-}(\partial_z^l f_N^K, f_N^K).
\end{equation*}
For the gain term, multiplying $\partial_z^l f_N^K \mathbf{1}_{\left\lbrace \partial_z^l f_N^K \leq 0\right\rbrace}$ and integrating on $ v \in \D $ and $z \in  I_z$ gives 
\begin{equation}\label{QR_plus}
\begin{split}
&\int_{I_z} \int_{\D} Q^{R,+} \left(\partial_z^l f_N^K, f_N^K \right)(t,v,z) \, \partial_z^l f_N^K \mathbf{1}_{\left\lbrace \partial_z^l f_N^K \leq 0\right\rbrace}(t,v,z) \,\rd v \,\pi(z) \,\rd z \\[2pt]
=& \int_{I_z}\int_{\D} Q^{R,+} \left(\partial_z^l f_N^{K,+} - \partial_z^l f_N^{K,-},  f_N^{K,+} - f_N^{K,-}\right)(t,v,z) \, \partial_z^l f_N^K \mathbf{1}_{\left\lbrace \partial_z^l f_N^K \leq 0\right\rbrace}(t,v,z) \,\rd v \,\pi(z) \,\rd z \\[2pt]
=& \int_{I_z} \int_{\D} \Big[ Q^{R,+}\left(\partial_z^l f_N^{K,+}, f_N^{K,+} - f_N^{K,-}\right) \\
&\qquad \qquad\qquad - Q^{R,+}\left(\partial_z^l f_N^{K,-}, f_N^{K,+} - f_N^{K,-}\right)\Big] \, \partial_z^l f_N^K \mathbf{1}_{\left\lbrace \partial_z^l f_N^K \leq 0\right\rbrace} \,\rd v \,\pi(z) \,\rd z \\[2pt]
=& \int_{I_z}\int_{\D} \Big[ Q^{R,+}\left(\partial_z^l f_N^{K,+}, f_N^{K,+}\right) - Q^{R,+}\left(\partial_z^l f_N^{K,+}, f_N^{K,-}\right) \\[2pt]
&\qquad \qquad\qquad - Q^{R,+}(\partial_z^l f_N^{K,-},  f_N^{K,+}) + Q^{R,+}(\partial_z^l f_N^{K,-}, f_N^{K,-})\Big] \left( -\partial_z^l f_N^{K,-} \right) \,\rd v \,\pi(z) \,\rd z \\[2pt]
=& \int_{I_z} \int_{\D}\Big[ - Q^{R,+}\left(\partial_z^l f_N^{K,+}, f_N^{K,+}\right) + Q^{R,+}\left(\partial_z^l f_N^{K,+}, f_N^{K,-}\right) \\[2pt]
&\qquad \qquad\qquad + Q^{R,+}\left(\partial_z^l f_N^{K,-}, f_N^{K,+}\right) - Q^{R,+}\left(\partial_z^l f_N^{K,-}, f_N^{K,-}\right)\Big] \, \partial_z^l f_N^{K,-}  \, \rd v \,\pi(z) \,\rd z, 
\end{split}
\end{equation}
so the gain part can be further estimated as follows:
\begin{equation}
    \begin{split}
    &\int_{I_z} \int_{\D} Q^{R,+}\left(\partial_z^l f_N^{K} , f_N^{K})(t,v,z \right) \, \partial_z^l f_N^{K} \mathbf{1}_{\left\lbrace \partial_z^l f_N^{K} \leq 0\right\rbrace}(t,v,z) \,\rd v \,\pi(z) \,\rd z\\
    \leq & \int_{I_z} \int_{\D}\left[ Q^{R,+}(\partial_z^l f_N^{K,+},  f_N^{K,-}) + Q^{R,+}(\partial_z^l f_N^{K,-},  f_N^{K,+}) \right](t,v,z) \, \partial_z^l f_N^{K,-}(t,v,z)\,\rd v \,\pi(z) \,\rd z\\
    \leq &\ \| Q^{R,+}( f_N^{K,+},  f_N^{K,-})\|_{L_{v}^2 H^1_z} \| \partial_z^l f_N^{K,-}\|_{L_{v}^2 L_z^2} 
    + \| Q^{R,+}( f_N^{K,-},  f_N^{K,+})\|_{L_{v}^2 H^1_z} \| \partial_z^l f_N^{K,-}\|_{L_{v}^2 L_z^2} \\
    \leq & \  C^+_{R,L,I_z,d,p}(B) \, \left( \| f_N^{K,+} \|_{L_v^1 H_z^1} \| f_N^{K,-} \|^2_{L_v^2 H_z^1} + \| f_N^{K,-} \|_{L_v^1 H_z^1} \| f_N^{K,+} \|_{L_v^2 H_z^1} \| f_N^{K,-} \|_{L_v^2 H_z^1} \right) \\
    \leq & \  C_{0}(B) \, \left( \| f_N^{K} \|_{L_v^1 H_z^1} \| f_N^{K,-} \|^2_{L_v^2 H_z^1} + \| f_N^{K} \|_{L_v^2 H_z^1} \| f_N^{K,-} \|^2_{L_v^2 H_z^1} \right),
    \end{split}
\end{equation}
where the estimate for the gain term \eqref{QGLp1} is used.\\ 
For the loss term,
\begin{equation}\label{QR_neg_final}
\begin{split}
&- \int_{I_z} \int_{\D} Q^{R,-}(\partial_z^l f_N^{K}, f_N^{K})(t,v,z) \, \partial_z^l f_N^{K} \mathbf{1}_{\left\lbrace \partial_z^l f_N^{K} \leq 0 \right\rbrace}(t,v,z)\, \rd v \,\pi(z) \,\rd z \\[4pt]
\leq & \int_{I_z} \int_{\D}  \left( L^{R}[\partial_z^l f_N^{K,+}]f_N^{K,-} \right) \partial_z^l f_N^{K,-} + \left( L^{R}[\partial_z^l f_N^{K,-}] f_N^{K} \right) \partial_z^l f_N^{K,-}  \,\rd v \,\pi(z) \,\rd z\\[6pt]
\leq & C^-_{R,L,I_z,d,p}(B) \left( \|  f_N^{K,+} \|_{L_v^1 H_z^1}  \| f_N^{K,-} \|_{L_v^2 H_z^1} \| \partial_z^l f_N^{K,-} \|_{L_{v,z}^2} + \|  f_N^{K,-}  \|_{L_v^1 H_z^1}  \| f_N^{K} \|_{L_v^2 H_z^1} \| \partial_z^l f_N^{K,-} \|_{L_{v,z}^2} \right) \\[6pt]
\leq & C_{0}(B) \left( \| f_N^{K} \|_{L_v^1 H_z^1}
\| f_N^{K,-} \|^2_{L_v^2 H_z^1}  + \| f_N^{K}  \|_{L_v^2 H_z^1} \| f_N^{K,-} \|^2_{L_v^2 H_z^1} \right), 
\end{split}
\end{equation}
where the estimate for the loss term \eqref{QLLp} is used. 

To estimate another term in \eqref{QR_Z_split}, namely $Q^R(f_N^{K}, \partial_z^l f_N^{K})$, the calculation is almost the same as above and we omit the details here. 

Combining \eqref{gn_char}, \eqref{EN}, \eqref{QR_plus}--\eqref{QR_neg_final} and adding up $|l|=0,1$, we obtain, for any $t \in [0,t_0]$,
\begin{equation}
\begin{split}
	\frac{1}{2} \frac{\rd}{\rd t} \left\| f_N^{K,-} \right\|^2_{L_v^2 H_z^1} 
	& \leq C_{0}(B) \Big( \| f_N^{K}\|_{L_v^1 H_z^1}  + \| f_N^{K}\|_{L_v^2 H_z^1} \Big) \left\| f_N^{K,-} \right\|^2_{L_v^2 H_z^1}  \\[4pt]
        & \qquad  + \frac{C_{\mP}C_{1}(B)\, K^2_{1,1}}{N} \left\| f_N^{K,-} \right\|_{L_v^2 H_z^1} \, ,
	\end{split}
\end{equation} 
which further yields that
\begin{equation}
\begin{split}
		\frac{\rd}{\rd t} \left\| f_N^{K,-} \right\|_{L_v^2 H_z^1} 
		& \leq C_{0}(B) \left( \| f_N^{K}\|_{L_v^1 H_z^1} + \| f_N^{K} \|_{L_v^2 H_z^1} \right) 
		 \left\| f_N^{K,-}\right\|_{L_v^2 H_z^1} + \frac{C_{\mP}C_{1}(B)\, K^2_{1,1}(t_0)}{N}\\[4pt]
		& \leq C_{0}(B) \left( 2E_1^{f^0} + K_{0,1}(t_0)  \right) \left\| f_N^{K,-} \right\|_{L_v^2 H_z^1} + \frac{C_{\mP}C_{1}(B)\, K^2_{1,1}(t_0)}{N} \\[4pt]
		& =: \mathcal{C}_{K_{0,1}(t_0)} \left\| f_N^{K,-}\right\|_{L_v^2 H_z^1}
		+ \frac{\mathcal{C}_{K_{1,1}(t_0)}}{N}\,. 
	\end{split}
\end{equation} 
This completes the proof of \eqref{fe_minus} by further applying the Gr\"onwall's inequality.
\end{proof}

We remark that the estimate \eqref{fe_minus} will be valid for the numerical solution $f_N^K$ obtained in our main well-posedness, i.e., Theorem \ref{existencetheorem} below.

\subsection{Proof of the well-posedness of the numerical solution \texorpdfstring{$f_N^K$}{fNK} with uncertainties}
\label{subsec:wellposed}

\subsubsection{Local well-posedness of the numerical solution \texorpdfstring{$f_N^K$}{fNK} on a small time interval \texorpdfstring{$[t_0, t_0+\tau]$}{[t0+tau]]}}
\label{subsub:local}

We first present the local well-posedness of the numerical solution $f_N^K$ in a well-constructed space.

\begin{proposition}\label{localexistence}
    Let the truncation parameters $R$, $L$ satisfy \eqref{RL} and assume the collision kernel $B$ satisfies \eqref{kernel}--\eqref{cutoff} and \eqref{kinetic}--\eqref{Assump_B}. 
    If we further assume that the initial condition $f^0(v,z)$ in the continuous system \eqref{ABE} belongs to $L^{1}_{v}H^1_z \cap L^{2}_{v}H^1_z$ with the following quantities
	\begin{equation}
	E^{f^0}_{1}=\|f^0(\cdot,\cdot)\|_{L^{1}_{v}H^1_z},  \quad  D^{f^0}_{0,1}=\left\|f^{0}(\cdot,\cdot)\right\|_{L^{2}_{v}H^1_z},
	\end{equation}
	and that the numerical system \eqref{PFSNK} starts to evolve from a certain time $t_0$ with
	\begin{equation}\label{fN_IC}
	\|f_N^K (t_0,\cdot,\cdot)\|_{L^{1}_{v}H^1_z} \leq 2 E^{f^0}_{1}, \quad \|f_N^K(t_0,\cdot,\cdot)\|_{L^{2}_{v}H^1_z} \leq K_{0,1}.
	\end{equation}
    Then there exists a local time $ \tau$ such that \eqref{PFSNK} admits a unique solution $ f_{N}^K = f_{N}^K(t,v,z) \in L^{1}_{v}H^1_z \cap L^{2}_{v}H^1_z $ on $ [t_0,t_0+\tau]$. One can particularly choose
	\begin{equation} \label{tau}
		\tau=\frac{1}{2(\bar{C}_1 \bar{D}+\bar{C}_2\bar{E})}, \quad \text{with} \quad \bar{E} = 4E^{f^0}_{1},  \quad \bar{D} = 2K_{0,1},
	\end{equation}
	such that
	\begin{equation}
		\forall t\in [t_0,t_0+\tau], \quad  \|f_N^K (t,\cdot,\cdot)\|_{L^{1}_{v}H^1_z} \leq \bar{E}, \quad  \|f_N^K(t,\cdot,\cdot)\|_{L^{2}_{v}H^1_z} \leq \bar{D},
	\end{equation}
    where the constants $\bar{C}_1$, $\bar{C}_2$ depend on the collision kernel $B$, domain of random variable $I_z$, the truncation parameters $R$, $L$, and dimension $d$.
\end{proposition}
\begin{proof}	
    The main strategy is to use the Banach fixed point theorem. The proof is almost the same as our preceding work \cite[Proposition 4.3]{HQY21} and \cite[Proposition 4.5]{LQ2022}, except one needs to apply the new estimate of the collision operator as in Corollary \ref{sumQz}.

\end{proof}

\subsubsection{Well-posedness of the numerical solution \texorpdfstring{$f_N^K$}{fNK} on an arbitrarily prescribed time interval \texorpdfstring{$[0,T]$}{[0,T]}}
\label{subsub:main}

In this subsection, we are ready to extend our local well-posedness of the numerical solution $f_{N}^K$ obtained in Proposition \ref{localexistence} to any arbitrarily prescribed time.

\begin{theorem}\label{existencetheorem}
    Let the truncation parameters $R$, $L$ satisfy \eqref{RL} and assume the collision kernel $B$ satisfies \eqref{kernel}--\eqref{cutoff} and  \eqref{kinetic}--\eqref{Assump_B}.
    If the initial condition $f^{0}(v,z)$ in the continuous system \eqref{ABE} and the numerical approximation $f_N^{0,K}(v,z)$ in the semi-discretized system \eqref{PFSNK} satisfy the assumptions specified in \eqref{con(a)}--\eqref{con(d)}.
	
    Then there exist integers $N_0$ and $K_0$ depending on the prescribed final time $T$ and initial condition $f^0$, such that for all $ N>N_{0} $ and $ K>K_{0} $, the semi-discretized system \eqref{PFSNK} admits a unique solution $f_N^K = f_N^K(t,\cdot,\cdot) \in L^1_v H^1_z \cap H^1_v H^1_z $ on the time interval $ [0,T] $ with the following estimates:
	\begin{equation}\label{fNL1L2}
	\forall t\in [0,T], \quad \left\| f_{N}^K(t,\cdot,\cdot) \right\|_{L^1_v H^1_z} \leq 2E^{f^0}_1, \quad \left\| f_{N}^K(t,\cdot,\cdot) \right\|_{H^1_{v}H^1_z} \leq K_{1,1}(T),
	\end{equation}
    where $K_{1,1}(T) >0$ is a constant depending on the prescribed time $T$, collision kernel $B$, domain of random variable $ I_z $, truncation parameters $R,L$, dimension $d$, and initial quantity $D^{f^0}_{0,1}=\|f^0(\cdot,\cdot)\|_{L^2_v H^1_z}$.
\end{theorem}

\begin{proof}
The time iteration is applied to extend the local well-posedness result in Proposition \ref{localexistence} to any arbitrarily prescribed time interval $[0,T]$.

\textbf{Step 1:} Start with the initial time $t=0$. By using condition \eqref{con(c)}, we can select large enough $N_1$ and $K_0$ such that for any $N \geq N_1$ and $K \geq K_0$,
\begin{equation} \label{initial}
    \|f^{0,K}_N(\cdot,\cdot)\|_{L^1_v H_z^1} \leq 2 E^{f^0}_1.
\end{equation}
Moreover, we have $\|f^{0,K}_N(\cdot,\cdot)\|_{L^2_v H_z^1} \leq D^{f^0}_{0,1} \leq K_{0,1}(T)$ thanks to the condition \eqref{con(b)}, where $K_{0,1}(T)$ can be explicitly given by the prescribed time $T$, collision kernel $B$, domain of random variable $ I_z $, truncation parameters $R,L$, dimension $d$, and initial quantity $D^{f^0}_{0,1}$. 

Then, by applying local well-posedness Proposition~\ref{localexistence}, there exists a unique solution $f_N^K(t,\cdot,\cdot)\in L^1_v H^1_z \cap L^2_v H^1_z$ over the time interval $[0,\tau]$ with the following estimate,
\begin{equation}
\forall t\in [0,\tau], \quad \|f_N^K(t,\cdot,\cdot)\|_{L^1_v H^1_z} \leq 4E^{f^0}_1.
\end{equation}
Furthermore, by taking advantage of the boundedness in $L^1_v H^1_z$ and that $f_N^{0,K}(\cdot,\cdot) \in H^1_v H^1_z$ from \eqref{con(b)}, we can invoke Proposition~\ref{regularity} to find the $L^2_v H^1_z$- and $H^1_v H^1_z$- estimates of the numerical solution $f_N^K$ in the local time interval $[0, \tau]$,
\begin{equation}
\forall t \in [0, \tau], \quad \|f_N(t,\cdot,\cdot)\|_{L^2_v H^1_z}\leq K_{0,1}(\tau), \quad \|f_N(t,\cdot,\cdot)\|_{H^1_v H^1_z} \leq K_{1,1}(\tau),
\end{equation}
as well as the estimate of the negative part $f^{K,-}_{N}$ in $L^2_v H^1_z$ from Proposition \ref{Pro_fnk_neg} that
\begin{equation}\label{fN-}
\begin{split}
\forall t\in [0, \tau], \quad \left\| f_N^{K,-}(t,\cdot, \cdot) \right\|_{L^2_v H^1_z} \leq e^{\tau \,\mathcal{C}_{K_{0,1}}} \left( \| f_N^{0,K,-}(\cdot, \cdot) \|_{L^2_v H^1_z} + \frac{\mathcal{C}_{K_{1,1}}}{N} \right).
\end{split}
\end{equation}

Furthermore, noticing that $|f_N^K(t,v,z)| = 2f_N^{K,-}(t,v,z) + f_N^K(t,v,z)$, we have

\begin{equation}\label{f-L1H1}
	\begin{split}
		\|f_N^K(t,\cdot,\cdot)\|_{L^1_v H^1_z}= &\| \int_{\D} f_N^K(t,v,\cdot) \,\rd{v} \|_{H^1_z} \\[4pt]
        =& \| \int_{\D} 2 f_N^{K,-}(t,v,\cdot) \,\rd{v} \|_{H^1_z} + \| \int_{\D} f_N^{K}(t,v,\cdot)  \,\rd{v} \|_{H^1_z} \\[4pt]
		=&2 \|f_N^{K,-}(t,\cdot,\cdot)\|_{L^{1}_{v}H^1_z} + \| \int_{\D} f^{0,K}(v,\cdot)  \,\rd v \|_{H^1_z} \\[4pt]
		\leq& 2 C_{L,I_z,d} \|f_N^{K,-}(t,\cdot,\cdot)\|_{L^{2}_{v}H^1_z}+ E^{f^0}_{1},
	\end{split}
\end{equation}
where $C_{L,I_z,d} > 0$ is a constant depending on the truncation parameters $L$, domain of random variable $ I_z $, and dimension $d$, and the conservation property \eqref{mass_conse_PK} in Lemma \ref{lemma:conv} is also utilized in the third equality above. 

Therefore, observing the estimate \eqref{f-L1H1}, it implies that $\|f_N^K(t,\cdot,\cdot)\|_{L^1_v H_z^1}$ will be under control as long as we have a ``good" estimate for the negative part $\|f_N^{K,-}(t,\cdot,\cdot)\|_{L^2_v H_z^1}$.
In fact, thanks to the estimate \eqref{fN-}, we can simply choose a $N_2$ large enough such that, for all $N \geq N_2$ and any prescribed time $T$, we can define a quantity $\bar{K}$ to satisfy
\begin{equation}\label{N0}
\begin{split}
\bar{K}:= \e^{ \mathcal{C}_{K_{0,1}(T)} \, T} \left( \| \fe^{0,K,-}(\cdot, \cdot) \|_{L_v^2H_z^1} + \frac{\mathcal{C}_{K_{1,1}(T)}}{N} \right) \leq \frac{E^{f^0}_1}{2 C_{L,I_z,d}},
\end{split}
\end{equation}
where the inequality above is feasible to achieve thanks to the condition \eqref{con(d)} and always holds for any 
$t \leq T$, as the quantity $\bar{K}$ is an increasing function with respect to time.

Hence, by combining \eqref{fN-}-\eqref{f-L1H1}, we find that
\begin{equation} \label{fNL1L21}
\forall t\in [0,\tau], \quad \|f_N(t,\cdot,\cdot)\|_{L_v^1H_z^1} \leq 2E^{f^0}_1.
\end{equation}


By choosing $N_0 = \max\{N_1, N_2\}$, where $N_1$ is determined to satisfy \eqref{initial} and $N_2$ is determined to satisfy \eqref{N0}, we have found such integers $N_0$ and $K_0$, depending only on the prescribed final time $T$ and initial condition $f^0$ such that for all $N > N_0$ and $K > K_0$, the numerical system \eqref{PFSNK} admits a unique solution $f_N^K(t,\cdot,\cdot) \in L_v^1 H_z^1 \cap L_v^2H_z^1$ in the time interval $[0,\tau]$, which satisfies \eqref{fNL1L21}. 

\textbf{Step 2:} By letting $t=\tau$ as a new initial time, we need to check if the local well-posedness can be extended in the equally-long time interval $[\tau,2\tau]$. Recalling the Step 1, where we have already proved that
\begin{equation}\label{initial1}
\forall t\in [0,\tau], \quad f_N^K(t,\cdot,\cdot) \in L_v^1 H_z^1 \cap L_v^2H_z^1 \quad \text{and} \quad \|f_N^K(t,\cdot,\cdot)\|_{L_v^1 H_z^1}\leq 2E^{f^0}_1.
\end{equation}
Then by taking $k=0$ and $r=1$ in the Proposition~\ref{regularity}, we have the following estimate for the propagation of $\|f_N^K(t,\cdot,\cdot)\|_{L_v^2 H_z^1}$,
\begin{equation}\label{initial2}
\begin{split}
    \forall t \in [0,\tau],  \quad \|f_N^K(t,\cdot,\cdot)\|_{L_v^2 H_z^1} \leq K_{0,1}(\tau) \leq K_{0,1}(T)
\end{split}
\end{equation}
From the estimate \eqref{initial1} - \eqref{initial2}, we find that $\|f_N^K(\tau,\cdot,\cdot)\|_{L_v^1 H_z^1}$ and $\|f_N^K(\tau,\cdot,\cdot)\|_{L_v^2 H_z^1}$ satisfy the condition of the local well-posedness Proposition~\ref{localexistence}, which allows us to apply Proposition~\ref{localexistence} starting from $t=\tau$ and further obtain that 
there exists a unique solution $f_N^K(t,\cdot,\cdot) \in L_v^1 H_z^1 \cap L_v^2H_z^1$ on $[\tau, 2\tau]$ with
\begin{equation}
\forall t\in [\tau, 2\tau], \quad \|f_N(t,\cdot,\cdot)\|_{L_v^1 H_z^1}\leq 4E^{f^0}_1.
\end{equation}

Meanwhile, noticing the bounded property in $L_v^1 H_z^1$ above and the fact that $f_N^{0,K}(\cdot,\cdot) \in H_v^1 H_z^1$, we can invoke the Proposition~\ref{regularity} over the interval $[0,2\tau]$ to derive that
\begin{equation}
\forall t \in [0, 2\tau], \quad \|f_N^K(t,\cdot,\cdot)\|_{L_v^2 H_z^1}\leq K_{0,1}(2\tau), \quad \|f_N^K(t,\cdot,\cdot)\|_{H_v^1 H_z^1}\leq K_{1,1}(2\tau),
\end{equation}
and for any $t\in [0,2\tau]$, 
\begin{equation} 
\begin{split}
	\left\|f_N^{K,-}(t,\cdot,\cdot)\right\|_{L_v^2 H_z^1}  \leq \e^{ (2\tau) \mathcal{C}_{K_{0,1}(2\tau)}} \left( \| f_N^{0,K,-}(\cdot, \cdot) \|_{L_v^2 H_z^1} + \frac{\mathcal{C}_{K_{1,1}(2\tau)}}{N} \right) \leq \bar{K}
\end{split}
\end{equation}
i.e., the same choice of $N$ chosen above would still make
\begin{equation} 
\forall t\in [0,2\tau], \quad \|f_N^K(t,\cdot,\cdot)\|_{L_v^1 H_z^1} \leq 2E^{f^0}_1.
\end{equation}
That is, at time point $t=2\tau$, we are back to the situation \eqref{initial1} at $t=\tau$. In fact, we can generalize the same strategy to longer time interval $[0, n\tau]$ with the same choice $N \geq N_0$, $K \geq K_0$ and the quantity $\bar{K}$, in the sense that, for all $t\in [0, n\tau]$,
\begin{equation}
    \begin{split}
        \|f_N^K(t,\cdot,\cdot)\|_{L_v^2 H_z^1}\leq &  K_{0,1}(n\tau), \quad \|f_N^K(t,\cdot,\cdot)\|_{H_v^1 H_z^1}\leq K_{1,1}(n\tau),\\[4pt]
        \left\|f_N^{K,-}(t,\cdot,\cdot)\right\|_{L_v^2 H_z^1} \leq & \bar{K}, \qquad \qquad \|f_N^K(t,\cdot,\cdot)\|_{L_v^1 H_z^1} \leq 2E^{f^0}_1.
    \end{split}
\end{equation}

\textbf{Step 3:} Repeating the Step 2 for $n$ times until the time interval $[0,n\tau]$ covers the prescribed interval $[0,T]$, we can show that there exists a unique solution $f_N^K(t,\cdot,\cdot) \in L_v^1 H_z^1 \cap H_v^1 H_z^1$ on $[0,T]$ with
\begin{equation}
\forall t\in [0,T], \quad \|f_N^K(t,\cdot,\cdot)\|_{L_v^1 H_z^1} \leq 2E^{f^0}_1, \quad \left\| f_N^K(t,\cdot,\cdot) \right\|_{L_v^2 H_z^1} \leq K_{0,1}(T),
\end{equation}
where the estimate of $\left\| f_N^K(t,\cdot,\cdot) \right\|_{L_v^2 H_z^1}$ can be verified by taking $k=0$ and $r=1$ in Proposition~\ref{regularity} once again.
\end{proof}


\section{Convergence and stability of the semi-discretized system for Boltzmann equation with uncertainties}
\label{sec:conv}

In this section, we will prove the convergence of the proposed semi-discretized system \eqref{PFSNK} by taking advantage of the well-posedness and stability of the numerical solution $f_N^K$ established in the previous section.

For the continuous system \eqref{ABE} with a periodic, non-negative initial condition $f^0(v,z) $ in $ L_v^1 H_z^r \cap H_v^k H_z^r$ for some integer $k \geq 1$ and $r \geq 1$, there exists a unique global non-negative solution $f(t,v,z) \in H_v^k H_z^r$ with the estimate that $\|f(t,\cdot,\cdot)\|_{H_v^k H_z^r} \leq C_{k,r}(f^0)$ for $t \geq 0$ with $C_{k,r}(f^0) > 0$ depending on initial datum $f^0$, which can be shown by following the similar argument to handle the deterministic model as in \cite[Proposition 5.1]{FM11} coupled with our estimate of the collision operator including uncertainties in Corollary \ref{sumQz}. 

For the numerical system \eqref{PFSNK}, we consider the initial condition $f_N^{K,0}(v,z) = \mP_N^K f^{0}(v,z)$ with $f^0(v,z) $ in $ L_v^1 H_z^1 \cap H_v^1 H_z^1$, which then satisfies the four conditions \eqref{con(a)}--\eqref{con(d)}.
Then, by Theorem \ref{existencetheorem}, there exists a unique solution $f_N^K(t,\cdot,\cdot) \in L^{1}_{v}H^1_z \cap H^{k}_{v}H^1_z$ over the prescribed time interval $[0,T]$, where $\|f_N^K(t,\cdot,\cdot)\|_{L^{2}_{v}H^1_z}\leq K_{0,1}(T)$  and $\|f_N^K(t,\cdot,\cdot)\|_{H^{1}_{v}H^1_z}\leq K_{1,1}(T)$ for $t \in[0,T]$.

We can show the following main theorem of the spectral convergence, in the sense that the error function $e_{N}^K$ with uncertainties is defined as follows:
\begin{equation}
e_N^K (t,v,z) := \mathcal{P}_N^K f(t,v,z) - f_N^K (t,v,z).
\end{equation}
\begin{theorem}\label{spectralaccuracy}
Let the truncation parameters $R$, $L$ satisfy \eqref{RL} and assume the collision kernel $B$ satisfies \eqref{kernel}--\eqref{cutoff} and \eqref{kinetic}--\eqref{Assump_B}. 
If the initial condition $f^{0}(v,z)$ in the continuous system \eqref{ABE} belongs to $ L_v^1 H_z^r \cap H_v^k H_z^r$ for some integers $k \geq 1$ and $r \geq 1$ and the numerical approximation $f_N^{0,K}(v,z)$ in the semi-discretized system \eqref{PFSNK} satisfy the assumptions specified in \eqref{con(a)}--\eqref{con(d)}.

Choose $N_0, K_0$ to satisfy the conditions in Theorem \ref{existencetheorem}, then the numerical solution $f_N^K$ to the semi-discretized system \eqref{PFSNK} is convergent for all $N>N_0$ and $K>K_0$, and exhibits spectral accuracy in the sense that
\begin{equation}\label{final}
\forall t\in [0,T], \quad \left\| e_{N}^K(t,\cdot,\cdot) \right\|_{L_v^2 H_z^1} \leq C_{k,r}(T,f^0)\left(\frac{1}{N^{k}} + \frac{1}{K^{r-1}}\right),
\end{equation}
where $C_{k,r}(T,f^0) > 0$ depends on $k$, $r$, prescribed time $T$, initial condition $f^0(v,z)$ as well as the truncation parameters $R,L$, domain of random variable $ I_z $, dimension $d$, collision kernel $B$.
\end{theorem}

\begin{proof}
To obtain the differential equation satisfied by our defined error function $e_N^K$, we first take the projection operator $\mP_N^K$ to both hand sides of the continuous system \eqref{ABE}:
	\begin{equation}\label{PNKUP}
	\left\{
	\begin{array}{lr}  
	\partial_{t} \mathcal{P}_N^K f(t,v,z) = \mathcal{P}_N^K Q^{R}(f,f)(t,v,z),\\[4pt]
	f_{N}^{0,K}(v,z) = \mathcal{P}_N^K f^0(v,z), 
	\end{array}
	\right.
	\end{equation}
 
Then, by taking subtraction between \eqref{PFSNK} and \eqref{PNKUP}, we find that
	\begin{equation}\label{errorNK}
	\left\{
	\begin{aligned}
	&\partial_{t} e_{N}^K(t,v,z)= \mathcal{P}_{N}^K Q^{R}(f,f)(t,v,z) - \mathcal{P}_{N}^K Q^{R}(f_{N}^K,f_{N}^K)(t,v,z),\\
	& e_{N}^K (0,v,z) = 0.
	\end{aligned}
	\right.
	\end{equation}
	where the zero initial error comes with the fact that $f_N^K(0,v,z) = \mathcal{P}_N^K f^0(v,z) $ for the initial approximation.
 
Next, by multiplying both hand sides of \eqref{errorNK} by $e_N^K$, and integrating over $\D \times I_z$,
\begin{equation}\label{multi1}
\begin{split}
    \frac{1}{2} \frac{\rd}{\rd t} \left\| e_N^K \right\|^{2}_{L^{2}_{v,z}}  =& \int_{I_z} \int_{\D} \left[ \mathcal{P}_{N}^K Q^{R}(f,f) - \mathcal{P}_{N}^K Q^{R}(f_{N}^K,f_{N}^K)\right](t,v,z)  e_{N}^K (t,v,z) \,\rd v \,\pi(z) \,\rd z \\[4pt]
    \leq & \left\| \mathcal{P}_{N}^K \left( Q^{R}(f,f) - Q^{R}(f_{N}^K,f_{N}^K) \right) \right\|_{L^{2}_{v,z}} \left\|  e_{N}^K \right\|_{L^{2}_{v,z}}\\[4pt]
    \leq & C_{\mP} \left\| Q^{R}(f,f) - Q^{R}(f_{N}^K,f_{N}^K) \right\|_{L^{2}_{v,z}} \left\|  e_{N}^K \right\|_{L^{2}_{v,z}}\\[4pt]
    \leq & C_{\mP} \left\| Q^{R}(f,f) - Q^{R}(f_{N}^K,f_{N}^K) \right\|_{L_v^2 H_z^1} \left\| e_{N}^K \right\|_{L_v^2 H_z^1},
\end{split}
\end{equation}
furthermore, by taking $\partial_z$ to both hand sides of \eqref{errorNK}, multiplying by $\partial_z e_N^K$, and integrating over $\D \times I_z$, we obtain
\begin{equation}\label{multip1}
\begin{split}
    \frac{1}{2} \frac{\rd}{\rd t} \left\| \partial_z e_N^K \right\|^{2}_{L^{2}_{v,z}}  =& \int_{I_z} \int_{\D} \partial_z \left[ \mathcal{P}_{N}^K Q^{R}(f,f) - \mathcal{P}_{N}^K Q^{R}(f_{N}^K,f_{N}^K)\right](t,v,z)  \partial_ze_{N}^K (t,v,z) \,\rd v \, \pi(z) \,\rd z \\[4pt]
    \leq & \left\| \mathcal{P}_{N}^K \partial_z\left( Q^{R}(f,f) - Q^{R}(f_{N}^K,f_{N}^K) \right)\right\|_{L^{2}_{v,z}} \left\| \partial_z e_{N}^K\right\|_{L^{2}_{v,z}}\\[4pt]
    \leq & C_{\mP} \left\| \partial_z \left( Q^{R}(f,f) - Q^{R}(f_{N}^K,f_{N}^K) \right)\right\|_{L^{2}_{v,z}} \left\| e_{N}^K\right\|_{L_v^2 H_z^1}\\[4pt]
    \leq & C_{\mP} \left\| Q^{R}(f,f) - Q^{R}(f_{N}^K,f_{N}^K) \right\|_{L_v^2 H_z^1} \left\|  e_{N}^K \right\|_{L_v^2 H_z^1}.
\end{split}
\end{equation}
Combining \eqref{multi1}-\eqref{multip1} yields that, 
	\begin{equation}
	\frac{\rd }{\rd t} \left\| e_N^K \right\|_{L^2_vH^1_z} \leq C_{\mP} \left\| Q^{R}(f,f) - Q^{R}(f_{N}^K,f_{N}^K) \right\|_{L^2_vH^1_z}. 
	\end{equation}
Then, we estimate the right-hand side of the inequality above,
	\begin{equation}\label{de1}
	\begin{split}
	&\left\| Q^{R}(f,f) - Q^{R}(f_{N}^K,f_{N}^K)\right\|_{L^2_vH^1_z}\\[4pt]
	\leq & \left\|  Q^{R}(f-f_N^K,f)\right\|_{L^2_vH^1_z} + \left\|  Q^{R}(f_N^K,f-f_N^K) \right\|_{L^2_vH^1_z}\\[4pt]
    \leq & C_{0}(B) \|f - f_N^K \|_{L^1_vH^1_z} \|f\|_{L^2_vH^1_z} + C_{0}(B) \| f_N^K \|_{L^1_vH^1_z} \|f - f_N^K\|_{L^2_vH^1_z} \\[4pt]
	\leq & C_{0}(B) \left(\left\|f\right\|_{L^2_vH^1_z} + \left\| f_{N}^K\right\|_{L^2_vH^1_z} \right) \left\|  f - f_{N}^K \right\|_{L^2_vH^1_z} \\[4pt]
	\leq  &C_{0}(B)\left( C_{0,1}(f^0) + K_{0,1}(T)\right) \left\|  f - f_{N}^K \right\|_{L^2_vH^1_z},
	\end{split}
	\end{equation}
where the bi-linearity of collision operator $Q^R$ and estimate \eqref{Q_HkZ} are utilized in the first and second inequality, while the well-posedness Theorem \ref{existencetheorem} of the numerical solution $f_N^K$ and its associated estimate in $\|\cdot\|_{L_v^2 H_z^1}$ are applied in the last inequality above. Furthermore, we also have,
	\begin{equation}\label{Pf0}
	\begin{split}
	\left\|  f - f_{N}^K\right\|_{L^2_vH^1_z} 
	=& \left\|  f - \mathcal{P}_{N}^K f + \mathcal{P}_{N}^Kf - f_{N}^K\right\|_{L^2_vH^1_z} \\[4pt]
	\leq & \left\|  f - \mathcal{P}_{N}^Kf \right\|_{L^2_vH^1_z} +  \left\|  \mathcal{P}_{N}^K f - f_{N}^K \right\|_{L^2_vH^1_z} \\[4pt]
	\leq & \left\|  f - \mathcal{P}_{N}^K f  \right\|_{L^2_vH^1_z} + \left\|e_{N}^K \right\|_{L^2_vH^1_z}\\[4pt]
	\leq &  \mathcal{C}_{k,r,\mP} C_{k,r}(f^0) \left( \frac{1}{N^{k}} + \frac{1}{K^{r-1}} \right) + \left\|e_{N}^K \right\|_{L^2_vH^1_z},
	\end{split}
	\end{equation}
where we utilize the following estimate in the last inequality above,
\begin{equation}\label{Pf}
\begin{split}
    \left\|  f - \mathcal{P}_{N}^K f \right\|_{L^2_vH^1_z} =& \left\|  f - \mathcal{P}_{N}f + \mathcal{P}_{N}f- \mathcal{P}_{N}^K f \right\|_{L^2_vH^1_z}\\[4pt]
    =& \left\|  f - \mathcal{P}_{N}f \right\|_{L^2_vH^1_z} + \left\|   \mathcal{P}_{N}f - \mathcal{P}_{N}^K f \right\|_{L^2_vH^1_z}\\[4pt]
    \leq & \frac{ \mathcal{C}_{k,\mP_N} \|f\|_{H^{k}_{v}H^r_z}}{N^{k}} + C_{\mathcal{P}_N}\left\| f - \mathcal{P}^K f \right\|_{L^2_vH^1_z}\\[4pt]
    \leq & \frac{ \mathcal{C}_{k,\mP_N} \|f\|_{H^{k}_{v}H^r_z}}{N^{k}} + C_{\mathcal{P}_N} \frac{\mathcal{C}_{r,\mP^{K}} \|f\|_{H^{k}_{v}H^r_z}}{K^{r-1}}\\[4pt]
    \leq & \mathcal{C}_{k,r,\mP} \left( \frac{ \|f\|_{H^{k}_{v}H^r_z}}{N^{k}} + \frac{\|f\|_{H^{k}_{v}H^r_z}}{K^{r-1}} \right) \\[4pt]
    \leq & \mathcal{C}_{k,r,\mP} C_{k,r}(f^0) \left( \frac{1}{N^{k}} + \frac{1}{K^{r-1}} \right),
\end{split}
\end{equation}
 
Therefore, by substituting \eqref{Pf0}-\eqref{Pf} in \eqref{de1}, we have
\begin{equation*}
\begin{split}
     \frac{\rd}{\rd t} \left\| e_N^K (t,\cdot,\cdot) \right\|_{L^2_vH^1_z} \leq & \ C_{\mP} \, \mathcal{C}_{k,r,\mP} \, C_{0}(B)\left( C_{0,1}(T) + K_{0,1}(T)\right)\\[4pt]
     & \qquad \qquad \qquad  \left[ \left\|e_{N}^K(t,\cdot,\cdot) \right\|_{L^2_vH^1_z} + \mathcal{C}_{k,r,\mP} \, C_{k,r}(f^0) \left( \frac{1}{N^{k}} + \frac{1}{K^{r-1}} \right) \right]\\
     =: & \ C_{k,r,1}(T,f^0)  \left\|e_{N}^K(t,\cdot,\cdot) \right\|_{L^2_vH^1_z} + C_{k,r,2}(T,f^0)  \left( \frac{1}{N^{k}} + \frac{1}{K^{r-1}} \right),
\end{split}
\end{equation*}
which, by the Gr\"onwall's inequality, implies that, for all $t\in [0,T]$,
\begin{equation}
    \left\| e_N^K (t,\cdot,\cdot) \right\|_{L^2_vH^1_z} \leq \e^{C_{k,r,1}(T,f^0) \, T} \left[\left\|e_{N}^K (0,\cdot,\cdot) \right\|_{L^2_vH^1_z} + \frac{C_{k,r,2}(T,f^0)}{C_{k,r,1}(T,f^0)} \left( \frac{1}{N^{k}} + \frac{1}{K^{r-1}} \right) \right].
\end{equation}
Considering the zero initial error $e_N^K(0,v,z)\equiv 0$ in \eqref{errorNK}, we can finally obtain the result in \eqref{final} by denoting $C_{k,r}(T,f^0):=\e^{C_{k,r,1}(T,f^0) \, T} \frac{C_{k,r,2}(T,f^0)}{C_{k,r,1}(T,f^0)}$.

\end{proof}


\section{Conclusion}
\label{sec:con}

In this paper, we demonstrate the convergence with spectral accuracy of the semi-discretized numerical system for the Boltzmann equation with uncertainties. Herein, deterministic and random variables are simultaneously discretized using the Fourier-Galerkin spectral method and the gPC-based SG method, respectively.
In particular, we employ the energy estimate to prove the well-posedness of the numerical solution obtained by the semi-discretized system in the high-order Sobolev space that encompasses the velocity and random variables altogether. 
Beyond addressing uncertainties stemming from the collision kernel and initial conditions, future research can potentially involve uncertainty quantification originating from the boundary value for the inhomogeneous Boltzmann equation, and the investigation can even be extended to the more complicated multi-species system with random inputs.

\section*{Acknowledgement}
\label{sec:ack}
L.~Liu acknowledges the support by National Key R\&D Program of China (2021YFA1001200), Ministry of Science and Technology in China, Early Career Scheme (24301021) and General Research Fund (14303022 \& 14301423) funded by Research Grants Council of Hong Kong from 2021-2023.  
K.~Qi thanks Michael Herty and Jingwei Hu for their helpful discussions as well as the kind hospitality during his visit to RWTH Aachen University. 
	


\appendix

\section{Proof of Proposition \ref{pro_QLp}}\label{App-QLP}

\begin{proof}
(I) For the gain term $ Q^{R,+}(g,f) $, the proof is similar to the usual Boltzmann operator $ Q^+(g,f)$ on $\mathbb{R}^d$. However, due to the different computational domains in $v$ as well as the existence of random variable $z$, we need to restrict back to a bounded domain $\D \times I_z$. Therefore, we follow \cite[Theorem 2.1]{mouhot2004regularity} and \cite[Proposition 3.1]{HQY21} to give a complete proof of \eqref{QGLp1}.
Start with the definition of Lebesgue norm via duality,
\begin{equation*}
\left\|Q^{R,+}(g,f)(v,z)\right\|_{L^p_{v,z}} = \sup \left\lbrace \int_{\D \times I_z} Q^{R,+}(g,f)(v,z) \Psi(v,z) \,\rd v \,\rd z; \ \left\|\Psi(\cdot, \cdot)\right\|_{L^{p'}_{v,z}} \leq 1 \right\rbrace.
\end{equation*}
Applying the pre-post collisional change of variables, i.e., $ (v,v_{*},\sigma) \rightarrow \left(v',v_{*}', \widehat{(v-v_*)}\right) $ with $\widehat{(v-v_*)} := \frac{v-v_{*}}{|v-v_{*}|} $ and an unit Jacobian, we can obtain
\begin{small}
\begin{equation}\label{dual1}
\begin{split}
&\int_{\D\times I_z} Q^{R,+}(g,f)(v,z) \Psi(v,z) \,\rd v \,\rd z \\
=& \int_{\D\times I_z} \int_{\mathbb{R}^d} \left( \int_{\bS^{d-1}} \mathbf{1}_{|v-v_{*}|\leq R} \Phi(|v-v_*|) b(\sigma \cdot \widehat{(v-v_*)},z) \Psi(v',z) \,\rd \sigma \right) g(v_{*},z)f(v,z)\, \rd v_* \,\rd v \,\rd z\\
=&\int_{\D\times I_z} \int_{\mathcal{B}_{\sqrt{2}L+R}} \left( \int_{\bS^{d-1}} \mathbf{1}_{|v-v_{*}|\leq R} \Phi(|v-v_*|) b(\sigma \cdot \widehat{(v-v_*)},z) \Psi(v',z) \,\rd \sigma \right) g(v_{*},z)f(v,z)\, \rd v_* \,\rd v \,\rd z,
\end{split}
\end{equation}
\end{small}
where the second equality is obtained by noting that $|v_*|\leq |v|+|v-v_*|$ and that $v\in \D$ and $|v-v_*|\leq R$.

Next, we define the linear operator $ S $ as follows:
\begin{equation*}
\begin{split}
S\Psi(v,z) &= \int_{\bS^{d-1}} \mathbf{1}_{|v|\leq R} \Phi(|v|) b(\sigma \cdot \hat{v},z) \Psi\left(\frac{v+|v|\sigma}{2},z \right) \rd \sigma,
\end{split}
\end{equation*}
such that \eqref{dual1} can be written as, for $\tau_hf(v):=f(v-h)$,
\begin{small}
\begin{multline*}\label{dual2}
    \int_{I_z} \int_{\D} Q^{R,+}(g,f)(v,z) \Psi(v,z) \,\rd v \pi(z) \,\rd z \\
    = \int_{I_z} \int_{{\mathcal{B}_{\sqrt{2}L+R}} }g(v_{*},z) \left( \int_{\D}  f(v,z) (\tau_{v_{*}}S(\tau_{-v_{*}}\Psi))(v,z) \, \rd v  \right)\rd v_{*} \pi(z) \,\rd z.
\end{multline*} 
\end{small}
To derive the $ L^p_v$-estimate of operator $S$, we first study the operator $S$ in $L^1_v$ and $L^{\infty}_v$ norms. Denoting $v^+=\frac{v+|v|\sigma}{2}$, where $\left| v^+ \right| \leq |v|$, we find, for all $z \in I_z$,
\begin{equation*}
\|S\Psi(\cdot, z)\|_{L^{\infty}_{v}}\leq  \|b(\cdot, \cdot)\|_{L^1_{\sigma}(\mathbb{S}^{d-1})L^{\infty}_z} \|\mathbf{1}_{|v|\leq R} \Phi(|v|)\|_{L^{\infty}_v} \|\Psi(\cdot, z)\|_{L^{\infty}_{v}{(\mathcal{B}_{\sqrt{2}L} )}}.
\end{equation*}
Also
\begin{equation*}
\begin{split}
\|S\Psi(\cdot, z)\|_{L^1_{v}}&\leq  \|\mathbf{1}_{|v|\leq R} \Phi(|v|)\|_{L^{\infty}_v}  \int_{\D}\int_{\bS^{d-1}}  b(\sigma \cdot \hat{v},z) \left|\Psi(v^+,z) \right|\rd \sigma\,\rd{v} \\
&\leq \|\mathbf{1}_{|v|\leq R} \Phi(|v|)\|_{L^{\infty}_v}  \int_{\mathcal{B}_{\sqrt{2}L}} \int_{\bS^{d-1}}  b(\cos\theta, z) \left|\Psi\left(v^+, z\right) \right| \frac{2^{d-1}}{\cos^2\theta/2}\,\rd \sigma\,\rd{v^+} \\
& \leq C \|b(\cdot, \cdot)\|_{L^1_{\sigma}(\mathbb{S}^{d-1})L^{\infty}_z}\|\mathbf{1}_{|v|\leq R} \Phi(|v|)\|_{L^{\infty}_v} \|\Psi(\cdot, z)\|_{L^1_{v}(\mathcal{B}_{\sqrt{2}L})}.
\end{split}
\end{equation*}
Thanks to the Riesz-Thorin interpolation, we can further obtain the $L^p_v$-estimate,
\begin{equation*}
\|S\Psi(\cdot, z)\|_{L^p_{v}} \leq C_{R,L,I_z,d,p'}^+(B) \|\Psi(\cdot, z)\|_{L^p_{v}(\mathcal{B}_{\sqrt{2}L} )}, \quad  1 \leq p \leq \infty,
\end{equation*}
where $C_{R,L,I_z,d,p'}^+(B) = C^{1/p'} \|b(\cdot, \cdot)\|_{L^1_{\sigma}(\mathbb{S}^{d-1})L^{\infty}_z} \|\mathbf{1}_{|v|\leq R} \Phi(|v|)\|_{L^{\infty}_v}$. 

Then, noting $\text{Supp}(S\Psi) \subset \mathcal{B}_R \subset \D$,
\begin{equation*}
\begin{split}
&\left| \int_{\D \times I_z} Q^{R,+}(g,f)(v,z) \Psi(v,z) \,\rd v \pi(z) \,\rd z\right|\\
\leq & \int_{I_z} \int_{\mathcal{B}_{\sqrt{2}L+R}}  |g(v_{*},z)|  \left( \int_{\D}\left|f(v,z)\right| \left|(\tau_{v_{*}}S(\tau_{-v_{*}}\Psi))(v,z)\right| \,\rd v  \right) \rd v_{*} \pi(z) \,\rd z \\
\leq  &\int_{I_z} \int_{\mathcal{B}_{\sqrt{2}L+R}}\left|g(v_{*},z)\right| \left\|f(\cdot,z)\right\|_{L^{p}_v} \left\|\tau_{v_{*}}S(\tau_{-v_{*}}\Psi) (\cdot,z)\right\|_{L^{p'}_v} \,\rd v_{*} \pi(z) \,\rd z\\
\leq  &\int_{I_z} \int_{\mathcal{B}_{\sqrt{2}L+R}} \left|g(v_{*},z)\right| \left\|f(\cdot,z)\right\|_{L^{p}_v} \left\|S(\tau_{-v_{*}}\Psi)(\cdot,z)\right\|_{L^{p'}_v} \,\rd v_{*} \pi(z) \,\rd z\\
\leq  & C_{R,L,I_z,d,p}^+(B) \int_{I_z} \int_{\mathcal{B}_{\sqrt{2}L+R}} \left|g(v_{*},z)\right| \left\|f(\cdot,z)\right\|_{L^{p}_v} \left\|\tau_{-v_{*}}\Psi(\cdot,z)\right\|_{L^{p'}_v(\mathcal{B}_{\sqrt{2}L})} \,\rd v_{*} \pi(z) \,\rd z\\
\leq  & C_{R,L,I_z,d,p}^+(B) \int_{I_z}\int_{\mathcal{B}_{\sqrt{2}L+R}}\left|g(v_{*},z)\right| \left\|f(\cdot,z)\right\|_{L^{p}_v} \left\|\Psi(\cdot,z)\right\|_{L^{p'}_v(\mathcal{B}_{2\sqrt{2}L+R})} \rd v_{*} \pi(z) \,\rd z \\
= & C_{R,L,I_z,d,p}^+(B) \int_{I_z} \left\|g(\cdot,z)\right\|_{L^{1}_v(\mathcal{B}_{\sqrt{2}L+R})} \left\|f(\cdot,z)\right\|_{L^{p}_v} \left\|\Psi(\cdot,z)\right\|_{L^{p'}_v(\mathcal{B}_{2\sqrt{2}L+R})} \pi(z) \,\rd z \\
\leq & C_{R,L,I_z,d,p}^+(B) \int_{I_z} \left\|g(\cdot,z)\right\|_{L^{1}_v} \left\|f(\cdot,z)\right\|_{L^{p}_v} \left\|\Psi(\cdot,z)\right\|_{L^{p'}_v} \pi(z) \,\rd z\\
\leq & C_{R,L,I_z,d,p}^+(B) \left\|g(\cdot,\cdot)\right\|_{L^1_v L^{\infty}_z} \left\|f(\cdot,\cdot)\right\|_{L^p_{v,z}} \left\|\Psi(\cdot,\cdot)\right\|_{L^{p'}_{v,z}}\\
\text{or} \ \leq  & C_{R,L,I_z,d,p}^+(B) \left\|g(\cdot,\cdot)\right\|_{L^1_v L^{p}_z} \left\|f(\cdot,\cdot)\right\|_{L^p_v L^{\infty}_z} \left\|\Psi(\cdot,\cdot)\right\|_{L^{p'}_{v,z}},
\end{split}
\end{equation*}
which leads to the two estimates in \eqref{QGLp1}, respectively .\\

(II) For the loss term $Q^{R,-}$, we write it as $Q^{R,-}(g,f)(v,z)=L^R[g](v,z)f(v,z)$ with $L^R$ being a convolution operator defined by
    \begin{equation*}
    \begin{split}
    L^R[g](v,z) = &\|b(\cdot, z)\|_{L^1_{\sigma}(\mathbb{S}^{d-1})} \int_{\bR^d} \mathbf{1}_{|q|\leq R}\Phi(|q|) g(v-q) \,\rd{q} \\
    =& \|b(\cdot, z)\|_{L^1_{\sigma}(\mathbb{S}^{d-1})} \left(\mathbf{1}_{|v|\leq R}\Phi(|v|)\right)*g(v).
    \end{split}
    \end{equation*}
Then,
    \begin{equation*}
    \begin{split}
	\|Q^{R,-}(g,f)(\cdot,\cdot)\|_{L^p_{v,z}} \leq & \left \| L^R[g](\cdot,\cdot)\right\|_{L^{\infty}_{v,z}} \|f(\cdot,\cdot)\|_{L^p_{v,z}}\\
    = & \|b(\cdot, \cdot)\|_{L^1_{\sigma}(\mathbb{S}^{d-1})L^{\infty}_z} \left\| \left(\mathbf{1}_{|v|\leq R}\Phi(|v|)\right)*g(\cdot) \right\|_{L^{\infty}_{v}} \|f(\cdot,\cdot)\|_{L^p_{v,z}}\\
    \leq & \|b(\cdot, \cdot)\|_{L^1_{\sigma}(\mathbb{S}^{d-1})L^{\infty}_z} \left \| \mathbf{1}_{|v|\leq R}\Phi(|v|)\right\|_{L^{\infty}_{v}} \left\|g(\cdot,\cdot)\right\|_{L^{1}_v(\mathcal{B}_{\sqrt{2}L+R}) L^{\infty}_z} \|f(\cdot,\cdot)\|_{L^p_{v,z}}\\
    \leq & C \|b(\cdot, \cdot)\|_{L^1_{\sigma}(\mathbb{S}^{d-1})L^{\infty}_z}  \left\| \mathbf{1}_{|v|\leq R}\Phi(|v|)\right\|_{L^{\infty}_{v}} \|g(\cdot,\cdot)\|_{L^1_v L^{\infty}_z} \|f(\cdot,\cdot)\|_{L^p_{v,z}}\\
    \leq & C_{R,L,I_z,d,p}^-(B) \|g(\cdot,\cdot)\|_{L^1_v L^{\infty}_z}\|f(\cdot,\cdot)\|_{L^p_{v,z}},
    \end{split}
    \end{equation*}
    on the other hand, we can also estimate the loss term $Q^{R,-}$ in the following way:
    \begin{equation*}
     \begin{split}
    \|Q^{R,-}(g,f)(\cdot,\cdot)\|_{L^p_{v,z}} \leq & \left \| L^R[g](\cdot,\cdot)\right\|_{L^{\infty}_{v}L^{p}_z} \|f(\cdot,\cdot)\|_{L^{p}_{v} L^{\infty}_z}\\
    = & \|b(\cdot, \cdot)\|_{L^1_{\sigma}(\mathbb{S}^{d-1})L^{\infty}_z} \left\| \left(\mathbf{1}_{|v|\leq R}\Phi(|v|)\right)*g(\cdot)\right\|_{L^{\infty}_{v}} \|f(\cdot,\cdot)\|_{L^{p}_{v} L^{\infty}_z}\\
    \leq & \|b(\cdot, \cdot)\|_{L^1_{\sigma}(\mathbb{S}^{d-1})L^{\infty}_z} \left \| \mathbf{1}_{|v|\leq R}\Phi(|v|)\right\|_{L^{\infty}_{v}} \left\|g(\cdot,\cdot)\right\|_{L^{1}_{v}(\mathcal{B}_{\sqrt{2}L+R})L^{p}_z} \|f(\cdot,\cdot)\|_{L^{p}_{v} L^{\infty}_z}\\
    \leq & C \|b(\cdot, \cdot)\|_{L^1_{\sigma}(\mathbb{S}^{d-1})L^{\infty}_z} \left \| \mathbf{1}_{|v|\leq R}\Phi(|v|)\right\|_{L^{\infty}_{v}} \|g(\cdot,\cdot)\|_{L^{1}_{v}(\mathcal{B}_{\sqrt{2}L+R})L^{p}_z} \|f(\cdot,\cdot)\|_{L^{p}_{v} L^{\infty}_z}\\
    \leq & C_{R,L,I_z,d,p}^-(B) \|g(\cdot,\cdot)\|_{L^{1}_{v}L^{p}_z} \|f(\cdot,\cdot)\|_{L^{p}_{v} L^{\infty}_z}.
    \end{split}
    \end{equation*}

    Finally, considering the decomposition $Q^{R}(g,f) = Q^{R,+}(g,f) - Q^{R,-}(g,f) $, the estimate \eqref{QLp} can be derived by combining \eqref{QGLp1}-\eqref{QLLp} directly.
\end{proof}


\bibliographystyle{amsplain}
\bibliography{Qi_bibtex}
	
\end{document}